\magnification=\magstep1
\input amstex
\documentstyle{amsppt}
\catcode`\@=11 \loadmathfont{rsfs}
\def\mycal{\mathfont@\rsfs}
\csname rsfs \endcsname \catcode`\@=\active

\vsize=6.5in

\topmatter 
\title Constructing MASA{\text{\rm s}} \\ with prescribed properties \endtitle
\author  Sorin Popa \endauthor

\rightheadtext{Constructing MASAs}

\affil     {\it  University of California, Los Angeles} \endaffil

\address Math.Dept., UCLA, Los Angeles, CA 90095-1555\endaddress
\email  popa\@math.ucla.edu\endemail

\thanks Supported by NSF Grant DMS-1400208, a Simons Fellowship, and Chaire d'Excellence de la FSMP 2016\endthanks

\abstract  \ We consider an iterative procedure for constructing \ maximal \ abelian \newline $^*$-subalgebras (MASAs) satisfying 
prescribed properties in II$_1$ factors. This method  
pairs well with  the intertwining by bimodules technique and with properties of the MASA and of the ambient factor that can 
be described locally. We obtain such a local characterization for II$_1$ factors 
$M$ that have an {\it s-MASA}, $A\subset M$ (i.e., for which $A \vee JAJ$ is maximal abelian in $\Cal B(L^2M)$), and use this strategy    
to prove that any factor in this class has uncountably 
many non-intertwinable singular (respectively semiregular)  s-MASAs.  

\endabstract 

\endtopmatter

\document

\heading 0. Introduction \endheading

Given a separable II$_1$ factor $M$, one can construct a 
{\it maximal abelian $^*$-subalgebra} (abreviated hereafter as {\it MASA}) $A$ in $M$  as an  inductive 
limit of  finite partitions.  
This iterative procedure pairs well with properties of MASAs that can be characterized locally,     
allowing the construction of $A$ 
in a manner that makes  ``more and more'' of the desired properties be satisfied.

This  technique has been initiated in [P81a],  [P81d] where it was used to prove that any separable II$_1$ factor 
$M$ contains a  MASA $A\subset M$ whose {\it normalizer} $\Cal N_M(A):=\{u\in \Cal U(M)\mid 
uAu^*=A\}$ generates a factor ($A$ is {\it semiregular} in $M$; see [P81a]), as well as a  
MASA $A\subset M$ whose normalizer is trivial, i.e. $\Cal N_M(A)=\Cal U(A)$  ($A$ is {\it singular} in $M$; see [P81d]). 

In this paper we obtain more refined applications of this method, by combining it with two additional ingredients: the intertwining by bimodule 
technique ([P01], [P03]) and local 
properties of the ambient II$_1$ factor $M$, such as existence of non-trivial central sequences (i.e., 
property Gamma of [MvN43]) and {\it s-thin approximation}, a property that we introduce here and which will be defined shortly.  

Recall in this respect that if $Q, P$ are von Neumann subalgebras in a II$_1$ factor $M$, then 
we write $Q\prec_M P$ if there exists a Hilbert $Q-P$ sub-bimodule $\Cal H\subset L^2M$ such that $\text{\rm dim} \Cal H_P < \infty$. 
In certain cases (notably if $Q, P$ are MASAs)  
this condition is equivalent to the existence of a non-zero partial isometry $v\in M$ 
such that $v^*v\in Q$ and $vQv^*\subset P$. 

Our first result shows  that any separable II$_1$ factor $M$ 
contains an uncountable family of singular (respectively semiregular) 
MASAs $\{A_i\}_i$ such that $A_i \not\prec_M A_j$, $\forall i\neq j$, with $A$ containing non-trivial central sequences of $M$ whenever 
$M$ does.  This will in fact  follow  from the following stronger result.

\proclaim{0.1. Theorem} Let $M$ be a separable ${\text{\rm II}}_1$ factor $M$ and $N\subset M$ a subfactor with trivial relative commutant, 
$N'\cap M=\Bbb C$. Let  $P_n\subset M$  be a sequence of von Neumann subalgebras such that $N\not\prec_MP_n$, $\forall n$. 
Then $N$ contains a singular   
$($respectively semiregular$)$ maximal abelian $^*$-subalgebra $A$ of $M$ such that $A\not\prec_MP_n$, $\forall n$. 
Moreover, if $N\simeq R$, then one can take $A$ so that to satisfy $\Cal N_M(A)''=N$, and if $N$ contains non-trivial 
central sequences of $M$, then $A$ can be 
taken so that to contain non-trivial central sequences of $M$ as well. 
\endproclaim

We then consider the class of II$_1$ factors $M$ which have an {\it s-MASA}, 
i.e., a MASA $A\subset M$ such that the von Neumann algebra 
$A\vee JAJ\subset \Cal B(L^2M)$, generated by left and right multiplication 
by elements in $A$ on the Hilbert space $L^2M$, is a MASA in $\Cal B(L^2M)$.  We obtain a 
local characterization of factors in this class, by proving that $M$ has an s-MASA if and only if it 
satisfies the following approximation property, that we call {\it s-thin}: for  
any finite partition $\{p_i\}_i \subset M$, any finite set $F\subset M$ and any $\varepsilon >0$, 
there exist a partition $\{q_j\}_j\subset M$  refining $\{p_i\}_i$ 
and an element $\xi \in M$ such that  any $x\in F$ can be $\varepsilon$-approximated in the norm-$\| \ \|_2$ 
by linear combinations of elements of the form $q_j \xi q_k$.  
We show that factors with s-MASAs are closed under amplifications and inductive limits and 
combine their local characterization with the iterative procedure  
to prove the following:

\proclaim{0.2. Theorem}  If $M$ has an s-MASA, 
then there exist uncountably many non-intertwinable s-MASAs in $M$, 
which in addition can be chosen singular $($resp. semiregular$)$. 
\endproclaim

The typical example of s-MASAs in II$_1$ factors are the {\it Cartan} (or {\it regular}) MASAs, i.e., MASAs $A\subset M$  
for which $\Cal N_M(A)''=M$ (cf [FM77]). Any group measure space II$_1$ factor $M=L^\infty(X)\rtimes \Gamma$, 
obtained from a free ergodic 
measure preserving action $\Gamma \curvearrowright X$  
of a countable group $\Gamma$ on a probability measure space $(X, \mu)$,  
has $A=L^\infty(X)$ as a Cartan subalgebra, which is thus also an s-MASA. 
The above result shows that such factors necessarily have singular s-MASAs as well. Note that when $M$ is hyperfinite, 
this fact was already known since 
([D54], [Pu61]), where the first concrete exemples of singular s-MASAs were given. 

By [OP07], [PV11], [PV12], there are large classes of group measure space II$_1$ factors that have 
unique (up to unitary conjugacy) Cartan subalgebras ($=$ regular MASAs), while by Theorem 0.2 above, such a factor  
always has ``many'' non conjugate semiregular s-MASAs. 

There are by now several classes of II$_1$ factors known to have no Cartan subalgebras,  
obtained first by using free probability theory ([Vo96], then by using deformation-rigidity theory ([OP07], [CS11], [CSU11], [PV11], [PV12], [I12]). 
It is interesting to note that in each case when one could prove absence of Cartan MASAs by using free probability, 
the same techniques could be used to show absence of s-MASAs as well  
(notably for the free group factors $L(\Bbb F_n)$, cf. [G98]).

While there is much evidence that II$_1$ factors with s-MASAs but no Cartan subalgebras do exist, 
the problem of constructing such examples remains open. Another open problem is to find 
new proofs for the non-existence of s-MASAs in certain II$_1$ factors, such as the free group factors $L(\Bbb F_n)$. 
But perhaps the most ``urgent'' open problem in this direction is to find an intrinsic, local characterization of II$_1$ factors 
having Cartan subalgebras. Such an intrinsic characterization may lead to interesting applications in deformation-rigidity theory. 
It may  also allow to prove that the class of factors with Cartan MASAs is close to inductive limits, a permanence property 
that, as we mentioned above, factors with s-MASAs do have. We discuss these open problems and other related questions
in the last section of the paper. 

This work has been finalized while I was visiting RIMS and Kyoto University in September 2016. 
I am very grateful to Masaki Izumi and Narutaka Ozawa for the warm hospitality extended to me during my stay. 

\vskip.05in
{\it Added in the proof}. In their very recent paper {\it Thin} II$_1$ {\it factors with no Cartan subalgebras} (math.OA/1611.02138), 
Anna Krogager and Stefaan Vaes were  
able to construct a large class of II$_1$ factors that have s-MASAs but no Cartan subalgebras (in fact 
are even strongly solid) thus solving a problem stated above and in 5.1.2.   

\heading 1. Preliminaries \endheading

All finite von Neumann algebras that we consider in this paper will come with a fixed 
normal faithful trace state, denoted $\tau$, and they will always be assumed 
separable with respect to the Hilbert norm $\| \ \|_2$ implemented by $\tau$. If $M$ is a finite von Neumann algebra and $B\subset M$ is a von Neumann subalgebra, then 
$E_B$ denotes the unique $\tau$-preserving conditional expectation of $M$ onto $B$. If $B\subset M$ is merely a weakly closed $^*$-subalgebra of $M$
(so $1_B$ not necessarily equal to $1_M$), then 
we use the same notation $E_B$ for the unique trace preserving expectation of $1_BM1_B$ onto $B$ that preserves the 
trace state $\tau( \cdot )/\tau(p)$ on $1_BM1_B$. We denote by $\Cal U(M)$ the unitary group of a von Neumann algebra $M$.  
If $\Cal X$ is a Banach space and $S\subset \Cal X$ is a subset, then we denote 
by $(S)_1$ the set of elements in $S$ that have norm at most $1$. For all notations that are not specified in the text, 
we send the reader to the expository notes ([P06]), and for basics on von Neumann algebras to the classic book [D57] or the recent [AP17]. 

\vskip.05in

\noindent
{\bf 1.1. Perturbation of projections}. The following  result is well known, but we  
state it here in the specific form needed in this paper. 

\proclaim{1.1.1. Lemma}  Let $M$ be a finite von Neumann algebra, $B\subset M$ a diffuse von Neumann subalgebra and 
$e\in \Cal P(M)$. Then there exists a projection $f\in B$ of trace equal to $\tau(e)$ such that $\|f-e\|_2\leq 14\|e-E_B(e)\|_2 +\sqrt{13\|e-E_B(e)\|_2}.$
\endproclaim 
\noindent
{\it Proof}. By (Lemma 1.1 in [P81c]), if $p$ is the spectral projection 
of $E_B(e)$ corresponding to the interval $[1/2, 1]$, then $\|p-E_B(e)\|_2 \leq 13\|e-E_B(e)\|_2$.  

By using the Cauchy-Schwartz inequality, this implies 
$$
|\tau(p)-\tau(e)| = |\tau(p)-\tau(E_B(e))| \leq 13\|e-E_B(e)\|_2.  
$$
Thus, if we take $f\in \Cal P(B)$ so that $\tau(f)=\tau(e)$ and satisfying $f\leq p$ if $\tau(p)\geq \tau(e)$ and 
$f\geq p$ if $\tau(p)\leq \tau(e)$, then $|\tau(f)-\tau(p)|\leq 13\|e-E_B(e)\|_2$. Altogether, 
$$
\|f-e\|_2 \leq \|f-p\|_2+\|p-E_B(e)\|_2+\|E_B(e)-e\|_2 
$$
$$
\leq 14\|e-E_B(e)\|_2 +\sqrt{13\|e-E_B(e)\|_2}.
$$
\hfill$\square$

\noindent
{\bf 1.2. Embedding $L^\infty([0,1])$ in II$_1$ factors}. We will view a diffuse abelian von Neumann 
subalgebra $A$ of a separable II$_1$ factor $M$ as an embedding of $L^\infty([0, 1])\simeq A$ into $M$. 
Thus, if $L^\infty([0, 1])$ is represented as an inductive limit of finer and finer partitions (e.g., dyadic) generating the 
$\sigma$-algebra of Lebesgue measurable subsets of $[0,1]$, then 
such an embedding is determined by the corresponding increasing sequence of finite dimensional subalgebras $A_n \nearrow A$. 

As pointed out in [P13], any embedding $L^\infty([0, 1])\simeq A \subset M$  acts weak mixingly on $M\ominus A'\cap M$, and this entails 
the following 2-independence property:

\proclaim{1.2.1. Theorem [P13]}  Let $M$ be a finite von Neumann algebra and $B\subset M$ a diffuse von Neumann subalgebra. 
Given any finite set $F\subset M\ominus B\vee (B'\cap M)$, any  $n\geq 2$ and any $\varepsilon>0$, there exists a partition of $1$ 
with projections in $B$, $p_1, ..., p_n \in \Cal P(B)$ of trace $1/n$, such that $|\|p_i x p_i \|^2_2- \tau(p_i)^2\tau(x^*x)| \leq \varepsilon$, $\forall x\in F$. 
\endproclaim 

The fact that $L^\infty([0, 1])\simeq A\subset M$ is a MASA is an extremality condition for the embedding, which can be 
described locally as follows (see e.g. [P81a]): 

\proclaim{1.2.2. Lemma}  Let  $M$ be a separable 
finite von Neumann algebra and $B\subset M$ a von Neumann subalgebra. 
Let  $A_n\subset B$ be an increasing sequence of finite 
dimensional von Neumann subalgebras and denote $A=\overline{\cup_n A_n}^w$. 
Let $\{x_j\}_j \subset M$ be a countable set, $\| \ \|_2$-dense in the unit ball of $M$. 
Then $A$ 
is maximal abelian in $B$ if and only if $\lim_n \|E_{A_n'\cap M}(E_B(x_j))-E_{A_n}(x_j)\|_2 =0$, 
$\forall j\geq 1$. Moreover, if $A$ is maximal abelian in $B$, then  
$\lim_n \|E_{A_n'\cap (B \vee B'\cap M)}(x_j)-E_{A_n\vee B'\cap M}(x_j)\|_2 =0$, 
$\forall j\geq 1$. 
\endproclaim 
\noindent
{\it Proof.} Since $\{x_j\}_j$ dense in $(M)_1$ in the norm $\| \ \|_2$ implies $\{E_B(j_j)\}_j$ 
$\| \ \|_2$-dense in $(B)_1$, the first part amounts to $A$ being maximal abelian in $B$ iff $A'\cap B=A$. The last part follows from by combining the first part  with the fact that 
$\cap_n A_n'\cap (B \vee B'\cap M)=(A'\cap B)\vee (B'\cap M)$ while $A_n \vee (B'\cap M) \nearrow A\vee (B'\cap M)$. 
\hfill$\square$

\vskip .05in

Let us also mention a result  that's essentially contained in (Section A.1 of [P92]), 
but which we derive here from 1.2.1 and 1.2.2 above: 

\proclaim{1.2.3. Corollary}  Let $M$ be a separable finite von Neumann algebra and $B\subset M$ a diffuse von Neumann subalgebra. 
There exists a MASA $A$ in $B$  such that $A'\cap M= A \vee (B'\cap M)$. Moreover, if $B=N$ is a $\text{\rm II}_1$ 
factor, then one can take $A$ to be contained  in a hyperfinite $\text{\rm II}_1$ subfactor $A\subset R\subset N$  satisfying $R'\cap M=N'\cap M$. 
\endproclaim 
\noindent
{\it Proof.} Let $\{x_j\}_j \subset (M)_1$ be $\|  \ \|_2$-dense sequence in the unit ball of $M$. We construct recursively an increasing 
sequence of finite dimensional abelian von Neumann algebras $A_m\subset B$ such that 
$\|E_{A_m'\cap M}(x_j)-E_{A_m \vee B'\cap M}(x_j)\|_2 \leq 2^{-m}$ for all $1\leq j \leq m$. Assuming we have constructed these algebras 
up to $m=n$, we construct $A_{n+1}$ as follows. By Theorem 1.2.1, given any $\alpha>0$, 
there exists an abelian finite dimensional $^*$-subalgebra $A^0_{n+1}$ 
containing $A_n$ such that $\|E_{{A^0_{n+1}}'\cap M}(x_j)-E_{({A_{n+1}^0}'\cap B)\vee (B'\cap M)}(x_j)\|_2< \alpha$, $1\leq j \leq n+1$. 
Then by taking fist a MASA $A^0$ in $B$ that contains $A^0_{n+1}$ and then using Lemma 1.2.2, we find a finite dimensional 
abelian subalgebra $A^\alpha_{n+1}\subset A^0$ that contains $A^0_{n+1}$, such that 
$\|E_{({A^\alpha_{n+1}}'\cap B)\vee (B'\cap M)}(x_j)-E_{A^\alpha_{n+1}\vee B'\cap M}(x_j)\|_2< \alpha$ 
and $\|E_{{A^\alpha_{n+1}}'\cap M}(E_B(x_j))-E_{A^\alpha_{n+1}}(x_j)\|_2\leq 2^{-n-1}$, for all $1\leq j \leq n+1$. 
Taking $\alpha$ sufficiently small and letting $A_{n+1}=A^\alpha_{n+1}$, we get 
$$
\|E_{A_{n+1}'\cap M}(x_j)-E_{A_{n+1}\vee B'\cap M}(x_j)\|_2< 2\alpha \leq 2^{-n-1}, 
$$
while we still have $\|E_{A_{n+1}'\cap M}(E_B(x_j))-E_{A_{n+1}}(x_j)\|_2 \leq 2^{-n-1}$, for all $1\leq j \leq n+1$. 
But then $A=\overline{\cup_n A_n}^w\subset B$ clearly satisfies the required condition.  

In the case $B=N$ is a II$_1$ factor, then we construct the increasing sequence of abelian subalgebras $A_n$ above to also be dyadic 
(i.e., all its minimal projections be have the same trace, equal to some $2^{-k_n}$) and 
so that each $A_n$ be the diagonal of a matrix algebra $R_n$ 
with matrix units $\{e_{ij}^n\}_{i,j}$, such that each $e^n_{ij}$ is a sum of some $e^{n+1}_{kl}$. Assuming we have constructed 
$A_m \subset R_m=\text{\rm sp}\{e^m_{ij}\}_{i,j}$  for $1\leq m \leq n$, we construct it for $m=n+1$ by applying the first part  
to the inclusion $B_0=e^n_{11}Ne^n_{11}\subset e^n_{11}Me^n_{11}=M_0$ and the finite set $F_0=\{e^n_{1i}x_ke^n_{j1} \mid 1\leq k \leq n+1, 1\leq i,j \leq n\}$, 
with appropriate $\alpha$, to get a finite dimensional subalgebra $A^1_{n+1}\subset B_0$ 
such that  $\|E_{{A^1_{n+1}}'\cap M_0}(x)-E_{B_0\vee B_0'\cap M_0}(x)\|_2< \alpha$, $\forall x\in F_0$. Moreover, we can take $A^1_{n+1}$ 
to be dyadic. We then choose matrix units $\{e_{kl}\}_{k,l}\subset B_0$ such that $e_{kk}$ are the minimal projections of $A^1_{n+1}$ 
and then define $\{e^{n+1}_{st}\}_{s,t}=\{e^n_{i1}e_{kl}e_{1j}^n \mid k, l, 1\leq i,j \leq n \}$, $R_{n+1}=\text{\rm sp}\{e^{n+1}_{st}\}_{s, t}$. 

Thus, if we let $A=\overline{\cup_n A_n}^w$ as above and put $R=\overline{\cup_n R_n}^w$, then we still get the condition 
$A'\cap M = A \vee N'\cap M$, but also $R'\cap M=R'\cap (A'\cap M)=N'\cap M$. 

\hfill$\square$

\vskip .05in 

\noindent
{\bf 1.3. Intertwining subalgebras in factors}. We recall here some basic 
facts about the ``intertwining'' subordination relation between 
subalgebras in II$_1$ factors, from [P01], [P03].  We will follow the presentation [P05a] of this topic, which emphasized the ``intertwining space'' 
between subalgebras. 

Thus, if $M$ is a finite von Neumann algebra and $Q, P\subset M$ 
are weakly closed $^*$-subalgebras of $M$, then $\Cal I_M(Q, P)$ denotes the set of vectors $\xi \in L^2(1_QM1_P$ 
with the property that the Hilbert $Q-P$ bimodule $\overline{\text{\rm sp} Q\xi P}\subset L^2M$ has finite dimension as a right $P$-module. This space is clearly invariant 
to taking sums and to multiplication by $Q$ from the left and $P$ from the right.  
We call it the {\it intertwining $Q-P$ sub-bimodule} of $M$.

The space $\Cal I_M(Q, P)$ has left support $\leq 1_Q$ and right support $\leq 1_P$, 
it is invariant to multiplication from the left by $Q'\cap M$ and from the right 
by $P'\cap M$ and it is increasing in $P$ and decreasing in $Q$. Also,  
$\Cal I_M(Q, P)=\Cal I_M(Q_1, P_1)$ whenever 
$Q_1\subset Q$, $P\subset P_1$ have finite index in the sense of [PP84] 
(either  the ``probabilistic'' definition, or the existence of a finite orthonormal basis; see 1.2 in [P94] for 
the equivalence between these alternative definitions). 
Moreover, if $q\in Q$, $p\in P$ are projections that have central trace of support $1$ in $Q$, respectively $P$, then 
$\Cal I_M(qQq, pPp)=q\Cal I_M(Q, P)p$. 

We'll denote as usual  by $\langle M, P \rangle$ the basic construction algebra, defined as 
the commutant in $\Cal B(L^2(M1_P))$ of the algebra of right multiplication by elements in $P$. 
It is also equal to 
the von Neumann algebra generated by operators of the form $xe_Py^*$, with $x, y \in M1_P$, acting on $\xi \in \hat{M1_P}\subset L^2(M1_P)$ 
by $xe_Py^*(\xi)=xE_P(y^*\xi)$. 

Then  the projection $s_{Q,P}:=\vee \{s(\xi e_P \xi^*) \mid \xi \in \Cal I_M(Q, P)\}$ is equal to the support of the direct summand of 
$Q'\cap1_Q \langle M, P \rangle 1_Q$ generated by projections that are finite in $1_Q \langle M, P \rangle 1_Q$ (where we have used the notation $s(T)$ 
for the support projection of a positive operator $T$). 
Thus, if $\xi \in L^2M$, then $\xi \perp \Cal I_M(Q, P)$ iff $\xi e_P \xi^* s_{Q,P}=0$ and iff $\xi e_P \xi^*$ is orthogonal on any 
projection $q'\in Q'\cap \langle M, P \rangle$ with $q'\langle M, P \rangle q'$ finite. 

If $\Cal I_M(Q, P)\neq 0$, then we say that  {\it $Q$ can be intertwined into $P$ inside $M$}, 
and write $Q\prec_M P$. Theorem 2.1 in [P03] shows that this condition  is equivalent to the following: 
there exist projections $p\in P$, $q\in Q$, a unital
isomorphism $\psi: qQq \rightarrow pPp$ (not necessarily onto)
and a partial isometry $v\in
M$ such that $vv^*\in (qQq)'\cap qMq$, $v^*v \in \psi(qQq)'\cap pMp$,  
$xv = v\psi(x), \forall x\in qQq$, and $x \in qQq$, $xvv^*=0$, implies $x=0$. 
Justified by this 2nd characterization, one also uses the terminology 
{\it a corner of $Q$ can be embedded into $P$ inside $M$}  (cf. 2.4 in [P03]).

By (2.1 in [P03]), the relation $Q\prec_MP$ is also equivalent to the fact that the action Ad$\Cal U(Q)$ 
has a non-zero part that's ``compact relative to $P$''. This means by definition that 
the commutant of $Q$ in the semifinite von Neumann algebra $1_Q \langle M, P \rangle 1_Q$ 
contains non-zero finite projections or, equivalently, that the action Ad$\Cal U(Q) \curvearrowright L^2(1_Q\langle M, P \rangle 1_Q, Tr)$ 
has non-zero fixed points.  

By (1.3 in [P01]), $Q\prec_M P$ is also equivalent to the fact that $Q'\cap 1_Q \langle M, P \rangle 1_Q$ 
contains non-zero elements from the ideal $\Cal J(\langle M, P \rangle)$ of elements in 
$\langle M, P \rangle$ that are ``compact relative to $P$''. 

We will use the notation $Q\not\prec_M P$ when  the above conditions are not satisfied,  
i.e., when $\Cal I_M(Q, P)=0$. This means that the action Ad-action of $\Cal U(Q)$ on $L^2(1_Q \langle M, P \rangle 1_Q, Tr)$ is  
ergodic. With the terminology (2.9 in [P05b]), in the case $1_Q=1_M$ this amounts to Ad$\Cal U(Q)\curvearrowright M$ 
being {\it weak mixing relative to $P$}.  

We recall from (2.3 in [P03]) some useful necessary and sufficient criteria for the condition $Q\not\prec_M P$ to be satisfied. 

\proclaim{1.3.1. Theorem} Let $M$ be a finite
von Neumann algebra and $P,Q\subset M$ be weakly closed $^*$-subalgebras. For each $q\in \Cal P(Q)$, 
fix $\Cal U_q\subset \Cal U(qQq)$ a subgroup generating $qQq$ as a von Neumann algebra. The following conditions are 
equivalent: 
\vskip .05in
$(1)$ $Q\not\prec_M P$
\vskip .05in
$(2)$ There exists a total subset $X \subset M$ and 
a sequence $u_n \in \Cal U_1$ such that $\lim_n \|E_P(xu_ny)\|_2$ $ =0$, $\forall x, y  \in X$.
\vskip .05in 
$(3)$ Given any $q\in \Cal P(Q)$ 
there exists a sequence of unitary elements $u_n\in \Cal U_q$ such that $\lim_n \|E_P(xu_ny)\|_2$ $ =0$, $\forall x, y  \in M$. 
\vskip .05in 

Moreover, if $P$ is regular in $M$, then the above are also equivalent to: 

\vskip .05in
$(4)$ There exists a total subset $X \subset M$ and 
a sequence $u_n \in \Cal U_1$ such that $\lim_n \|E_P(xu_n)\|_2$ $ =0$, $\forall x  \in X$. 
\endproclaim

The proof of the above theorem in ([P03]) actually shows the following more general   
result, involving the intertwining space  (cf. [P05a]):

\proclaim{1.3.2. Theorem} With the same assumptions as in $1.3.1$, if 
$X \subset L^2(1_QM1_P)$, then the following conditions are equivalent:  
\vskip .05in
$(1)$ $X \perp \Cal I_M(Q, P)$. 
\vskip .05in
$(2)$ There exist $u_n \in \Cal U_1$ such that $\lim_n \|E_P(\xi^* u_n \xi)\|_1=0$, $\forall \xi \in X$. 
\vskip .05in
$(3)$ There exist $u_n \in \Cal U_1$ such that $\lim_n \|E_P(\eta^* u_n \xi)\|_1=0$, $\forall \xi \in X, \eta\in L^2M$. 
\endproclaim

\vskip.05in

\noindent
{\bf 1.3.3. Remarks}. $(a)$ Property $1.3.2 (2)$  above, for characterizing the orthogonal in $L^2M$ of the 
intertwining space $\Cal I_M(Q, P)$, can be traced  back to ([P81b]), where this type of condition appears 
in the case of subalgebras $Q=L(G_1), P=L(G_2)$ of $M=L(G)$, arising from subgroups $G_1, G_2\subset G$, 
as well as for general $M$ and $Q=P$ (as in 1.4 below).   

\vskip.05in

$(b)$  The relation $Q\prec_M P$ is a ``virtual'' subordination relation, 
in the sense that  it is ``insensitive to finite index perturbations'':  
if $Q$ or $P$ are replaced by subalgebras $Q_1\subset Q$, $P_1\subset P$ of finite index  
(in the sense of one of the definitions in ([PP84]), then we still have $Q_1\prec_MP_1$.  In particular, 
if $Q$ has a finite dimensional direct summand, then $Q\prec_MP$ for any $P\subset M$, and if there exist projections  
$p\in P$, $p'\in P'\cap M$ such that $pPpp'=pp'Mpp'$, then $Q\prec_M P$ for any $Q\subset M$. 
The relation $Q\prec_M P$ 
is also insensitive to localization to ``corners'' of the algebras involved, i.e., it is sufficient to be satisfied under 
non-zero projections of $Q$, $P$ (or of their commutants in $M$).  

\vskip.05in

$(c)$ Related to $(b)$ above, let us underline here that the notions of finite index (up to taking 
``corners'') for an inclusion of  finite von Neumann algebras, $P\subset M$   
considered in [PP84], generalizing the Jones index in the case of inclusions of factors [J83],  translates into 
the relation $M\prec_M P$.  More precisely, this last condition means that there exist projections 
$p\in P$, $p'\in P'\cap M$ such that $pPpp' = P_0 \subset M_0=pp'Mpp'$ has finite index, either in the sense that there exists a 
finite orthonormal basis of $M_0$ over $P_0$ ([PP84]) or that $E_{P_0}(x)\geq c x$, $\forall x\in (M_0)_+$, for some $c>0$ (see A.1 in [V07]). 

In turn, the opposite relation $M\not\prec_M P$ translates into the fact that $P$ has {\it uniform infinite index} in $M$
and it amounts to $\Cal U(pp'Mpp')$ containing sequences of elements that are ``more and more'' perpendicular to 
$pPpp'$, for any $p\in \Cal P(P)$, $p'\in \Cal P(P'\cap M)$.  This type of condition characterizing infinite index can be traced back to 
(2.2 in [PP84]). 

\vskip.05in
$(d)$ Since it is determined by its behavior on corners, the subordination relation $\prec_M$ is not transitive in general. 
For instance, if we take $Q=pMp+\Bbb C(1-p)$, $P=\Bbb Cp+ \Bbb C(1-p)$ then we have $M\prec_M Q$, $Q\prec_M P$, 
but $M\not\prec_M P$. For this same reason, requiring $Q\prec_MP$ and $P\prec_MQ$, 
does not define a ``reasonable'' equivalence relation $\sim_M$ between subalgebras of $M$ 
(e.g., the previous example would show that $M\sim \Bbb C$). However, for MASAs of $M$, we have the following (cf. A.1 in [P01]): 

\proclaim{1.3.4. Theorem} Let $M$ be a finite von Neumann algebra and $A, B\subset M$ be MASAs in $M$. Then 
$A\prec_M B$ if and only if $B\prec_MA$ and if and only if there exists a non-zero partial isometry $v\in M$ such that $vv^*\in B$, 
$v^*v\in A$ and $vAv^*=Bvv^*$. 
\endproclaim

\noindent
{\bf 1.4. Normalizing subalgebra}. If $M$ is a finite von Neumann algebra and $B\subset M$ is 
a von Neumann subalgebra, then we denote $\Cal N_M(B):=\{u\in \Cal U(M) \mid uBu^*=B\}$, the {\it normalizer of $B$ in $M$}. 
The von Neumann algebra it generates, $\Cal N_M(B)''$, is called the {\it normalizing von Neumann algebra of $B$ in $M$}. 

A von Neumann subalgebra $B$ is {\it singular in $M$} if any automorphism Ad$(u)$ implemented by some $u\in \Cal N_M(B)$ is 
inner, i.e., it is of the form Ad$(v)$ for some $v\in \Cal U(B)$. This is the same as requiring that $\Cal N_M(B)=\Cal U(B)\Cal U(B'\cap M)$. 
If in turn $\Cal N_M(B)''=M$, then we say that $B$ is {\it regular in $M$}. 

This terminology has been introduced in [D54], 
in the case $B=A\subset M$ is a maximal abelian $^*$-subalgebra (MASA) in $M$.  
Note that for a MASA $A\subset M$, being singular  means that $\Cal N_M(A)''=A$, or equivalently 
$\Cal N_M(A)=\Cal U(A)$, i.e., the normalizer of $A$ in $M$ acts trivially on $A$. 
A regular MASA will be called a {\it Cartan subalgebra} 
(or {\it Cartan MASA}) in $M$. We will also consider MASAs 
$A\subset M$ for which  $\Cal N_M(A)''$ is a factor (equivalently, $\Cal N_M(A)$ acts  ergodically on $A$), 
which will be called {\it semi-regular} (cf. [D54]).

More generally, recall from ([P97] and 1.4 in [P01]) that if $B\subset M$ is a von Neumann subalgebra then $q\Cal N_M(B)$ 
denotes the set of all $x\in M$ with the property that there exists $x_1, ..., x_n \in M$ such that $Bx \subset \Sigma_i x_i B$ and 
$xB \subset \Sigma_i Bx_i$. The space $q\Cal N_M(B)$ is a $^*$-subalgebra and we see that, by (Lemma 1.4.2 in [P01]), 
one has $q\Cal N_M(B)=\Cal I_M(B, B)\cap \Cal I_M(B, B)^*\cap M$, with the weak closure being a von Neumann 
subalgebra of $M$. 

Note that if $B=A$ is a MASA in $M$, 
then by (1.3 in [P01]) we have $q\Cal N(A) =\text{\rm sp}\Cal N_M(A) = \Cal I_M(A, A)\cap M$ and this space is $\| \ \|_2$-dense in 
$\Cal I_M(A, A)$. Thus, the normalizing von Neumann algebra of $A$ 
satisfies $\Cal N_M(A)''=\overline{q\Cal N_M(A)}^w$ and the orthogonal of this space in $L^2M$ 
coincides with $\Cal I_M(A, A)^\perp$.  So Theorem 1.3.2 entails the following criterion for estimating the 
size of the normalizer of $A$ in $M$:

\proclaim{1.4.1. Corollary} Let $M$ be a finite von Neumann algebra, $A\subset M$  a MASA and $N=\Cal N_M(A)''$ 
its normalizing von Neumann algebra. The following conditions are equivalent for an element $\xi \in L^2M$: 
\vskip .05in
$(1)$ $\xi \perp N$. 
\vskip .05in
$(2)$ $\exists \{u_n\}_n\subset \Cal U(A)$ such that $\lim_n \|E_A(\xi^* u_n\xi)\|_1=0$.
\vskip .05in
$(3)$ For any $n\geq 2$ and any $\varepsilon >0$, there exists an abelian von Neumann subalgebra $A_0\subset A$ generated by 
$n$ projections of equal trace such that $\|E_A(\xi^* y\xi)\|_1 \leq \varepsilon$, $\forall y\in \{A_0\ominus \Bbb C1 \mid \|y\|\leq 1\}$. 

\endproclaim
\noindent
{\it Proof.} By applying 1.3.2 to the case $Q=P=A$, and taking into account that for a MASA $A\subset M$ one has 
$\Cal I_M(A, A)^\perp = \Cal N_M(A)^\perp$, it follows that $(1) \Leftrightarrow (2)$. Then by taking $\Cal U_1\subset \Cal U(A)$ 
to be a subgroup satisfying $\Cal U_1''=A$ and $\tau(u)=0$,  $u^2=1$, $\forall u\in  \Cal U_1\setminus \{1\}$, 
we get $(2) \Leftrightarrow (3)$.  
\hfill $\square$

\vskip.05in
Finally, let us note that a singular MASA $A\subset M$ means an embedding 
of the diffuse abelian von Neumann algebra $L^\infty([0, 1]) \simeq A \subset M$ 
so that the Ad-action of its unitary group on $M$ is weak mixing relative to $A$, 
a Cartan MASA is an embedding so that this action is compact relative to $A$, while a semi-regular MASA 
is an embedding having a large relative compact part.

\heading 2. Constructing MASAs with control of intertwiners  \endheading

Results in [P81a], [P81d] show that 
any separable II$_1$ factor $M$ has semi-regular and singular MASAs. 
The proof consists in constructing an embedding of $L^\infty([0, 1])$  $\simeq A \subset M$ as an 
inductive limit of dyadic partitions $A_n \nearrow A$ 
that become ``more and more extremal in $M$'' (resulting into $A$ being a MASA), while also 
controlling  the normalizer of $A$, making it become singular (in [P81d]),  respectively semi-regular (in [P81a]). 

For $A$ to become singular, one needs $A_n$ to ``become more and more relative weak mixing''. For it to 
become semi-regular, it is sufficient to build $A_n$ so that fixed matrix units having $A_n$ as diagonal are in the normalizer 
(i.e., in the relative compact part), at each step $n$. 

We will show below how one can use much more 
of  the intertwining by bimodules criteria within such iterative procedure, 
allowing us to construct embeddings $L^\infty([0, 1])\simeq A \subset M$ 
so that to be weak mixing relative to a given countable family of subalgebras of $M$. One can in fact   
even control such relative weak mixingness when $M$ is embedded into larger II$_1$ factors, thus leading to super-rigidity type properties for $A$. 
Moreover, we will do the construction so that to also take  into account local properties of the ambient factor $M$, 
such as existence of central sequences (in this section), and s-thin approximation (in the next section).

\proclaim{2.1. Theorem} Let $N$ be a separable $\text{\rm II}_1$ factor and 
$N \hookrightarrow M_n$ be embeddings of $N$ into separable $\text{\rm II}_1$ factors 
such that $N'\cap M_n$ is of type $\text{\rm I}$, $\forall n$. Let also  
$P_n \subset M_n$ be von Neumann subalgebras. 

\vskip.05in
\noindent
$1^\circ$ There exists a MASA $A\subset N$ such that for each $n$ one has: 

\vskip.05in
\noindent
$(a)$ $\Cal N_{M_n}(A)=\Cal U(A \vee N'\cap M_n)$; 

\noindent
$(b)$ $\Cal I_{M_n}(A, P_n)^\perp = \Cal I_{M_n}(N, P_n)^\perp$; 

\noindent
$(c)$ $M_n'\cap A^\omega$ is non-trivial whenever $M_n'\cap N^\omega$ is non-trivial. 

\vskip.05in
\noindent
In particular, $A$ is singular in $N$ and if $N'\cap M_n=\Bbb C1$, then $A$ is a singular MASA in $M_n$ which  
contains non-trivial central sequences of $M_n$ whenever $N$ does.  

\vskip.05in
\noindent
$2^\circ$ There exists a semiregular MASA $A\subset N$ such that for each $n$ one has: 

\vskip.05in
\noindent
$(a)$ $A'\cap M_n=A \vee N'\cap M_n$; 

\noindent
$(b)$ $\Cal N_{M_n}(A)''\subset N\vee N'\cap M_n$;  

\noindent
$(c)$ $\Cal I_{M_n}(A, P_n)^\perp = \Cal I_{M_n}(N, P_n)^\perp$. 

\noindent
$(d)$ $M_n'\cap A^\omega$ is non-trivial whenever $M_n'\cap N^\omega$ is non-trivial. 

\vskip.05in
\noindent
Moreover, if $N\simeq R$ then one can take $A\subset N$ such that 
$\Cal N_{M_n}(A)''=N\vee N'\cap M_n$. 

\endproclaim 
\noindent
{\it Proof}. For each $M_n$ choose a sequence $\{x^n_k\}_k\subset (M_n)_1$ that's $\| \ \|_2$-dense in $(M_n)_1$ 
and a sequence $\{\xi^n_k\}_k \subset L^2(M_n)\ominus \Cal I_{M_n}(N, P_n)$ that's $\| \ \|_2$-dense in 
$L^2(M_n)\ominus \Cal I_{M_n}(N, P_n)$. 
 
Let also $\{e_m\}_m \subset \{ e\in \Cal P(N) \mid \tau(e)\leq 1/2\}$ be a $\| \ \|_2$-dense sequence. 

To prove $1^\circ$, we construct recursively 
a sequence of finite dimensional abelian von Neumann subalgebras $A_m\subset N$ together with projections 
$f_m\in \Cal P(A_m)$, with $\tau(f_m)=\tau(e_m)$, and unitary elements $v_m \in \Cal U(A_mf_m)$, $w_m, u_m \in \Cal U(A_m)$,  
satisfying the following properties for all $1\leq i,j,k \leq m$:
$$
\|f_m - e_m\|_2 \leq 13 \|e_m - E_{A_{m-1}'\cap N}(e_m)\|_2 \tag 2.1.1
$$
$$
\|E_{A_m'\cap M_k}({x^k_i}^*v_mx^k_j)(1-f_m)\|_2 \leq 2^{-m},  \tag 2.1.2
$$
$$
\|E_{A_m'\cap M_k}(x^k_j)-E_{A_m \vee N'\cap M_k}(x^k_j)\|_2 \leq 2^{-m},  \tag 2.1.3
$$
$$
\|E_{P_k}({x^k_i}^* w_m \xi^k_j)\|_2 \leq 2^{-m},  \tag 2.1.4
$$
$$
\|[x^k_i, u_m]\|_2 \leq 2^{-m}, \| E_{A_{m-1}}(u_m)\|_2 \leq 2^{-m}.  \tag 2.1.5
$$

Assume we have constructed $(A_m, f_m, v_m, w_m, u_m)$ satisfying these properties for $m=1, 2, ..., n$. By applying Lemma 1.1.1 to $B=A_n'\cap N$ and $e=e_{n+1}$, 
it follows that there exists $f_{n+1}\in A_n'\cap N$ such that  
$\|f_{n+1}-e_{n+1}\|_2\leq 13\|e_{n+1}-E_{A_n'\cap N}(e_{n+1})\|_2$ and $\tau(f_{n+1}=\tau(e_{n+1})$. By Corollary 1.2.3, there exists a MASA $B_0\subset (1-f_{n+1})N(1-f_{n+1})$ 
satisfying the property

$$
B_0'\cap (1-f_{n+1})M_k(1-f_{n+1})=B_0\vee 
(N'\cap M_k) (1-f_{n+1}), \forall 1\leq k \leq n+1.$$ 

Since $(A_nf_{n+1})'\cap f_{n+1}Nf_{n+1}$ is type II$_1$ and $B_0\vee (N'\cap M_k)(1-f_{n+1})$ are of type I, 
for each $k$ we have $(A_nf_{n+1})'\cap f_{n+1}Nf_{n+1}\not\prec_{M_k} B\vee (N'\cap M_k)(1-f_{n+1})$. Thus, there exists $v_{n+1} 
\in \Cal U((A_nf_{n+1})'\cap f_{n+1}Nf_{n+1})$ with the property that 
$$
\|E_{B_0'\cap M_k}((1-f_{n+1}){x^k_i}^*v_{n+1}x^k_j(1-f_{n+1}))\|_2 < 2^{-n-1}, 1\leq i,j,k \leq n+1. \tag 2.1.6 
$$
Moreover, we may clearly assume $v_{n+1}$ has finite spectrum. We then take a refinement $A^0_{n+1}$ of $A_n$ in $N$ that contains $f_{n+1}$,  
such that $A^0_{n+1}f_{n+1}$ contains $v_{n+1}$, while $A^0_{n+1}(1-f_{n+1})$ ``approximates'' $B_0$ well enough (in the sense of Lemma 1.2.1) 
so that,  due to $(2.1.6)$ and its strict inequality, we still have   
$$
\|E_{{A^0_{n+1}}'\cap M_k}((1-f_{n+1}){x^k_i}^*v_{n+1}x^k_j(1-f_{n+1}))\|_2 < 2^{-n-1}, 1\leq i,j,k \leq n+1.  \tag 2.1.7
$$

On the other hand, by Corollary 1.2.3, there exists a finite dimensional abelian von Neumann subalgebra $A^1_{n+1}$ in $N$ 
that contains $A^0_{n+1}$ and satisfies 
$$
\|E_{{A^1_{n+1}} '\cap M_k}(x^k_j)-E_{A^1_{n+1} \vee N'\cap M_k}(x^k_j)\|_2 \leq 2^{-n-1}, 1\leq i,j,k \leq n+1.  \tag 2.1.8
$$

Now, since ${A^1_{n+1}}'\cap N$ has finite index in $N$, by Section 1.3 we have $\Cal I_{M_k}({A^1_{n+1}}'\cap N, P_k)=\Cal I_{M_k}(N, P_k)$, 
and thus $\xi^k_j \perp \Cal I_{M_k}(({A^1_{n+1}}'\cap N, P_k)$, for all $1\leq j, k \leq n+1$. Thus, by Thereom 1.3.2,  there exists 
a unitary element $w_{n+1}\in ({A^1_{n+1}}'\cap N)$ such that
$$
\|E_{P_k}({x^k_i}^*w_{n+1}\xi^k_j)\|_2 < 2^{-n-1}, 1\leq i,j,k \leq n+1.  \tag 2.1.9
$$
Also, we may clearly assume $w_{n+1}$ has finite spectrum. 
We take $A^2_{n+1}\subset N$ to be the 
finite dimensional abelian von Neumann algebra generated by $A^1_{n+1}$ 
and $w_{n+1}$. Due to  $(2.1.7)$, $(2.1.8)$ and $(2.1.9)$, $A^2_{n+1}$ satisfies conditions $(2.1.1)-(2.1.4)$.  

Finally, 
by using the fact that $M_n'\cap N^\omega \neq \Bbb C$ implies $M_n'\cap N^\omega$ diffuse (because $M_n$ is a factor), it follows that for any 
$\alpha>0$ there exists a projection $p\in N$ of trace $1/2$ such that $\|[x, p]\|_2  < \alpha$ and $|\tau(px)-\tau(p)\tau(x)|<\alpha$  
for all $x\in \{x^k_i\mid 1\leq i, k \leq n+1\}\cup (A^2_{n+1})_1$. 
By taking $\alpha$ sufficiently small and using Lemma 1.1.1, it follows that there exists $u_{n+1}\in \Cal U({A^2_{n+1}}'\cap N)$ 
sufficiently close to $1-2p$ so that we have $\|[x^k_i, u_{n+1}]\|_2 \leq 2^{-n-1}$, for all $1\leq i,k \leq n+1$. 
Thus, if we define $A_{n+1}=A^2_{n+1}\vee \{u_{n+1}\}''$, then conditions $(2.1.1)-(2.1.5)$ are satisfied for $m=n+1$.  

Define $A=\overline{\cup_n A_n}^w$. Condition $(2.1.3)$ clearly implies $A'\cap M_k = A \vee N'\cap M_k$, $\forall k$, 
while condition $(2.1.4)$ implies $\Cal I_{M_n}(A, P_n)^\perp = \Cal I_{M_n}(N, P_n)^\perp$, $\forall k$. 
 
Let $u\in \Cal N_{M_m}(A)$. If Ad$(u)$ acts non-trivially on $A$, then there exists a non-zero projection $e\in A$ 
of trace $\leq 1/2$ such that $u^*eu \leq 1-e$. 
Let $n_0$ be large enough so that $2^{-n_0} < \|e\|_2/60$. Since $\{e_n\}_n$ is $\| \ \|_2$-dense in the set of projections of $N$ 
of trace $\leq 1/2$, there exists $n\geq n_0$ such that $\|e_n - e\|_2 < \|e\|_2/60$ and  such that there exists $j,k\leq n$ with 
$\|x^k_j-eu\|_2 < \|e\|_2/60$. With $f_n, v_n \in A_n\subset A$ as given by the construction, $(2.1.1)$ implies 
$$
\|f_n-e_n\|_2 \leq 13\|e_n-E_A(e_n)\|_2 = 13\|e_n-e +E_A(e)-E_A(e_n)\|_2
$$
$$
\leq 26\|e_n-e\|_2 \leq 26\|e\|_2/60. 
$$
Since $u^*ev_nu\in A\subset A_n'\cap M_m$, we would then get the estimates 
$$
\|ef_n\|_2=\|ev_n\|_2=\|E_{A_n'\cap M_m}(u^*ev_nu)(1-e)\|_2 \tag 2.1.10
$$
$$
\leq \|E_{A_n'\cap M_m}(u^*v_nu)(1-f_n)\|_2 + \|f_n-e_n\|_2 + \|e_n-e\|_2 
$$
$$
\leq  \|E_{A_n'\cap M_m}({x^k_j}^*v_nx^k_j)(1-f_n)\|_2 + 29\|e\|_2/60 \leq 30\|e\|_2/60,   
$$
where for the very last inequality we used $(2.1.1)$. But since $\|e_n-e\|_2 < \|e\|_2/60$, we also have 
$$
\|f_ne-e\|_2\leq \|f_n-e\|_2\leq \|f_n-e_n\|_2+\|e_n-e\|_2\leq 27 \|e\|_2/60, 
$$
which together with $(2.1.10)$ implies that  
$$
30\|e\|_2/60 \geq \|ef_n\|_2\geq \|e\|_2-27\|e\|_2/60=33\|e\|_2/60, 
$$
a contradiction. This shows that $\Cal N_{M_k}(A)=\Cal U(A'\cap M_k)=\Cal U(A \vee N'\cap M_k)$, $\forall k$, finishing the 
proof that $A$ satisfies all the conditions in part $1^\circ$ of the theorem. 

To prove $2^\circ$, note first that by Corollary 1.2.3 there exists a hyperfine II$_1$ subfactor $R\subset N$ 
such that $R'\cap M_k = N'\cap M_k$, $\forall k$. It is then sufficient to construct a Cartan subalgebra $A$ of $R$ such that $A'\cap M_k = A\vee N'\cap M_k$,  
$\Cal N_{M_k}(A)''=R \vee N'\cap M_k$, $\Cal I_{M_n}(A, P_n)^\perp = \Cal I_{M_n}(N, P_n)^\perp$, $\forall k$. In other words, it is 
sufficient to prove the last part of $2^\circ$, where one assumes $N\simeq R$ and want to construct a Cartan MASA $A\subset N$ 
whose normalizing algebra is exactly $N\vee N'\cap M_k$. 

To this end, we construct recursively 
a sequence of commuting dyadic matrix subalgebras 
$R_m\subset N$ (i.e., $R_m\simeq M_{2^{k_m}\times 2^{k_m}}(\Bbb C)$, for some $k_m\geq 1$), 
with diagonal subalgebras $D_m\subset R_m$, such that if we denote  
$N_m=\vee_{k=1}^m R_k$, $A_m=\vee_{k=1}^m D_k$, 
there exist a projection $f_m$ of trace $1/2$ in $D_m$ 
and unitary elements $v_m \in \Cal U(D_mf_m)$, $w_m \in \Cal U(D_m)$, 
so that if we denote $y^k_i=x^k_i-E_{N\vee N'\cap M_k}(x^k_i)$, then the following properties 
are satisfied for $1\leq i,j,k \leq m$:  
$$
\|E_{A_m'\cap M_k}({y^k_i}^*v_my^k_j)(1-f_m)\|_2 \leq 1/10; \tag 2.1.11
$$
$$
\|f_m(y^k_i)(1-f_m)\|_2\geq 2\|y^k_i\|_2/5; \tag 2.1.12
$$
$$
\|E_{A_m'\cap M_k}(x^k_j)-E_{A_m \vee N'\cap M_k}(x^k_j)\|_2 \leq 2^{-m};  \tag 2.1.13
$$
$$
\|E_{P_k}({x^k_i}^* w_m \xi^k_j)\|_2 \leq 2^{-m};  \tag 2.1.14
$$
$$
\|E_{N_m}(x^k_i)-E_N(x^k_i)\|_2 \leq 2^{-m}. \tag 2.1.15
$$

Assuming we have constructed these objects up to $m=n$, we construct them for $m=n+1$ as follows. 

Noticing that the finite set $F=\{y^k_j \mid 1\leq j,k \leq n+1\}$ is perpendicular 
to $N_n \vee N'\cap M_k$ (which is equal to the commutant in $M_k$ of $N_n'\cap N$), by Lemma 1.2.1 we can first pick 
a projection $f$ of trace $1/2$ in $N_n'\cap N$ such that $f$ is almost $2$-independent to 
$F$. In particular, we can choose $f$ so that  for all $1\leq j, k \leq n+1$ we have
$$
\|fy^k_j(1-f)\|_2\geq 2\|y^k_j\|_2/5. \tag 2.1.16
$$ 

By Corollary 1.2.3,  there exists a MASA $B_0$ in $(1-f)(N_n'\cap N)(1-f)$ such that 
$B_0'\cap (1-f)M_k(1-f)=B_0\vee (N_n\vee N'\cap M_k)(1-f)$, $\forall 1\leq k \leq n+1$.  
Since $f(N_n'\cap N)f$ is type II$_1$ and $B_0\vee (N_n \vee N'\cap M_k)(1-f)$ is type I, 
the former cannot be intertwined into the latter inside $M_k$, so by Theorem 1.3.2 there exists a unitary element $v\in f(N_n'\cap N)f$ such that 
for all $1\leq i,j,k \leq n+1$ we have
$$
\|E_{B_0'\cap (1-f)M_k(1-f)}((1-f){x^k_i}^*vx^k_j(1-f))\|_2 < 1/10. \tag 2.1.17
$$

Moreover, we may choose $v$ so that to belong to a dyadic finite dimensional abelian subalgebra $B_1^1\subset f(N_n'\cap N)f$.  
Also, by approximating $B_0$ sufficiently well 
with a dyadic finite dimensional subalgebra $B_1^0\subset B_0$,  we will still have for all $1\leq i,j,k \leq n+1$ the estimates 
$$
\|E_{{B_1^0}'\cap (1-f)M_k(1-f)}((1-f){x^k_i}^*vx^k_j(1-f))\|_2 < 1/10. \tag 2.1.18
$$

We now take $B_1\subset (B_1\vee N_n)'\cap N$ to be a dyadic finite dimensional abelian subalgebra containing $B_1^1f+B_1^0(1-f)$. 
Since $\xi^k_j \perp \Cal I_{M_k}(N, P_k) = \Cal I_{M_k}((B_1\vee N_n)'\cap N, P_k)$, there exists a unitary element $w\in  B_1'\cap N$ such that 
$$
\|E_{P_k}({\xi^k_i}^* w \xi^k_j)\|_2 \leq 2^{-n-1}, 1\leq i,j,k \leq n+1,   \tag 2.1.19
$$

We may clearly also assume $w$ lies in a dyadic finite dimensional abelian subalgebra $B_2\subset N_n'\cap N$ that contains $B_1$.  
Moreover, by using Corollary 1.2.3 again, we may also assume $B_2$ is so that  $A_{n+1}^0=A_{n}\vee B_2$ satisfies 

$$
\|E_{{A^0_{n+1}}'\cap M_k}(x^k_j)-E_{A^0_{n+1} \vee N'\cap M_k}(x^k_j)\|_2 \leq 2^{-n-1}. \tag 2.1.20
$$ 

Take now $R_{n+1}^0\subset N_n'\cap N$ 
to be a (dyadic) finite dimensional factor having $B_2$ as a diagonal algebra. Finally, since $N\simeq R$, there exists a dyadic finite dimensional 
factor $R_{n+1}\subset N_n'\cap N$ that contains $R_{n+1}^0$, such that if we define $N_{n+1}=N_n \vee R_{n+1}$ then 
$$
\|E_{N_{n+1}}(x^k_i)-E_N(x^k_i)\|_2 \leq 2^{-n-1}, 1\leq i,k \leq n+1. \tag 2.1.21
$$
Thus, if we take $D_{n+1}$ to be a diagonal of $R_{n+1}$ that contains $B_2$ and denote $A_{n+1}=A_n \vee D_{n+1}$, 
$N_{n+1}=N_n \vee R_{n+1}$, $f_{n+1}=f$, $v_{n+1}=v$, $w_{n+1}=w$, then $(2.1.16)-(2.1.21)$ insure that conditions $(2.1.11)-(2.1.15)$ 
are satisfied for $n+1$. 

Let now $R=\vee_k R_k=\overline{\cup_n N_n}^w$, $A=\vee_k D_k = \overline{\cup_n A_n}^w$. Condition $(2.1.15)$ clearly implies that $R=N$, 
while condition $(2.1.13)$ implies $A'\cap M_k=A\vee N'\cap M_k$ and $(2.1.14)$ implies $\Cal I_{M_n}(A, P_n)^\perp = \Cal I_{M_n}(N, P_n)^\perp$, $\forall k$. 

By construction, we have that $A$ is Cartan in $N$, so that $\Cal N_{M_k}(A)''$ contains $N\vee N'\cap M_k$. If this inclusion is strict for some $k$, 
then by the factoriality of $N$ there must exist an automorphism $\theta$ of $A$ implemented by a unitary 
$u\in \Cal N_{M_k}(A)$ such that $\theta \circ \text{\rm Ad}(v)$ acts freely on $A$, $\forall v\in \Cal N_N(A)$.  Thus, $E_{N\vee N'\cap M_k}(u)=0$. 

Let $x^k_j\in M_k$ be so that $\|x^k_j-u\|_2 \leq 1/20$. This implies that 
$\|y^k_j-u\|_2 \leq 1/20$ and that for each $n\geq k$ we have $\|{y^k_j}^*v_ny^k_j-u^*v_nu\|_2 \leq 1/10$,  
while by $(2.1.12)$ we also have 
$$
\|f_nu(1-f_n)\|_2\geq \|f_ny^k_j(1-f_n)\|_2- 1/10
\geq 2 \|y^k_j\|_2/5 - 1/10
$$  
Thus, since $u^*v_nu\in A$, by $(2.1.11)$ we get  
$$
1/10 \geq \|E_{A \vee N'\cap M_k}({y^k_i}^*v_ny^k_j)(1-f_n)\|_2
$$
$$
\geq \|(1-f_n)u^*v_nu(1-f_n)\|_2-1/10 = \|f_nu(1-f_n)\|_2-1/10 
$$
$$
\geq 2 \|y^k_j\|_2/5 - 1/5\geq 2/5-1/20-1/5=3/20, 
$$
which is a contradiction. 
\hfill$\square$

\vskip.05in

Recall that in [P81d] one proves existence of singular MASAs not only in II$_1$ factors, but also in II$_\infty$ and 
III$_\lambda$ factors, for $0< \lambda <1$. This result is obtained as a consequence of stronger statement about 
MASAs in a (separable) II$_1$ factor $M$, showing that given any group of automorphisms 
$\Cal G$ of the associated II$_\infty$ factor $M^\infty=M\overline{\otimes} \Cal B(\ell^2\Bbb N)$ 
such that $\Cal G/\text{\rm Int}(M^\infty)$ is countable, there exists a $\Cal G$-{\it singular} MASA 
$A\subset M$, i.e., a maximal abelian subalgebra of $M$ with the property that if $\theta\in \Cal G$  
satisfies $\theta(a)\subset A$, for all $a$ in a ``corner'' $Ap$ of $A$, then $\theta$ 
acts as the identity on $Ap$. Let us note here that this type of result is in fact covered by the above general theorem: 

\proclaim{2.2. Corollary}  Let $M$ be a $\text{\rm II}_1$ factor with a sequence 
of von Neumann subalgebras of uniform infinite index $P_n \subset M$ $($in the sense of $1.3.3. (c))$. Let also $\Cal G\subset \text{\rm Aut}(M^\infty)$ 
be a subgroup of automorphisms of $M^\infty$ that contains $\text{\rm Int}(M^\infty)$ and is so that $\Cal G/\text{\rm Int}(M)$ is countable. 
Then there exists a $\Cal G$-singular MASA $A\subset M$ such that $A\not\prec_MP_n$, $\forall n$. In particular, $M$ contains  
uncountably many non-intertwinable $\Cal G$-singular MASAs, which in addition can 
be taken to contain non-trivial central sequences of $M$ whenever $M$ has property Gamma. 
\endproclaim 
\noindent
{\it Proof}. Let $\theta_n\in \Cal G$ be a sequence of automorphisms of $M^\infty$ such that $\Cal G=\cup_n \theta_n\circ \text{\rm Int}(M^\infty)$. 
For each $n$, let $s_n$ be so that $Tr\circ \theta_n = s_n Tr$ and denote $t_n=1+s_n$. 
Let $M_n=M^{t_n}$ and $f_n \in M_n$ a projection of trace $\tau(f_n)=1/(1+s_n)$. Thus, $f_nM_nf_n\simeq M$ and 
we can embed $M$ into $M_n$ as the subfactor $\{x \oplus \theta_n(x) \mid x\in M\simeq f_nM_nf_n\}$. 

By part $1^\circ$ of Theorem 2.1, there exists a singular MASA $A\subset M$ such that $A'\cap M_n = Af_n + A(1-f_n)$ and $\Cal N_{M_n}(A)=
\Cal U(Af_n + A(1-f_n))$. 
It is immediate to see that this means $A$ is $\Cal G$-singular in $M$. Moreover, by 2.1.1$^\circ$ we can take $A$ so that to also 
satisfy  $A\not\prec_{M_n} P_n$ (when $P_n$ is viewed as a subalgebra of $M_n$). This of course implies $A\not\prec_M P_n$ as well. 

To prove the last part, let $\Cal F$ be a maximal family of $\Cal G$-singular MASAs of $M$ such that $A, B$ are not intertwinable 
for any $A\neq B$ in $\Cal F$. Assume $\Cal F$ is countable and note that 
$M\not\prec_{M_n} A$, $\forall A\in \Cal F$. By applying Theorem 2.1 to $N=M\subset M_n$ (with $M_n$ as above), $\forall n$, and $\{P_n\}_n=\Cal F$,  
it follows that there exists a singular 
MASA $C\subset M$ such that $C\not\prec_M A$, $\forall A\in \Cal F$, contradicting the maximality of $\Cal F$. 

\hfill$\square$

\proclaim{2.3. Corollary}  Any separable $\text{\rm II}_1$ factor contains an uncountable family of mutually non-conjugate  
semi-regular MASAs, which in addition can be chosen to contain non-trivial central sequences if the ambient factor has 
property Gamma. 
\endproclaim 
\noindent
{\it Proof}. The same argument as in the above proof applies, using 2.1.2$^\circ$ instead of 2.1.1$^\circ$. 
\hfill$\square$

It has been shown in (2.5 of [P81b]) that if $P\subset M$ is a  
von Neumann subalgebra and $u\in \Cal U(M)$ satisfies the property that for all $n\geq 1$ and all $\varepsilon >0$ there exists   
a subalgebra $A_0\simeq L(\Bbb Z/n \Bbb Z) \subset Q$ with $u^* A_0 u \perp_\varepsilon P$, then $u \perp \Cal N_M(P)$. 
Along these lines, one can now deduce the following stronger result (generalizing 1.4.1 as well):  

\proclaim{2.4. Corollary}  Let $M$ be a finite von Neumann algebra. Let $Q, P \subset M$ be diffuse von Neumann 
subalgebras and $\xi \in L^2M$. The following conditions are equivalent:

$1^\circ$  $\xi \perp \Cal I_M(Q, P)$. 

$2^\circ$  For any $n\geq 1$ and any $\varepsilon >0$ 
there exists $u \in \Cal U(Q)$ such that $u^n=1$, $\tau(u^k)=0$, $1\leq k <n$ and $\|E_P(\xi^*u^k\xi)\|_1\leq \varepsilon$, $1\leq k<n$. 

$3^\circ$ For any $n\geq 1$ and any $\varepsilon > 0$, there 
exists an $n$-dimensional abelian von Neumann subalgebra $A_0\subset Q$ such that $\tau(p)=1/n$ for any minimal projection in $A_0$ 
and $\|E_P(\xi^*x\xi)\|_1\leq \varepsilon$ for all $x\in (A_0)_1$ with $\tau(x)=0$. 

\endproclaim 
\noindent
{\it Proof}. We clearly have $3^\circ \Leftrightarrow 2^\circ$ and by Theorem 1.3.2 we have $2^\circ \Rightarrow 1^\circ$. 

To prove $1^\circ \Rightarrow 3^\circ$, let $\xi\perp \Cal I_M(Q, P)^\perp$. Apply first Theorem 2.1 
to get a MASA $A\subset Q$ with the property that $\Cal I_M(A, P)^\perp=\Cal I_M(Q, P)^\perp$. 
Representing $A$ as the von Neumann algebra of the countable group $\Bbb Z/2\Bbb Z^{\oplus \infty}$, we get a unitary group $\Cal U_1\subset \Cal U(A)$ 
such that $\Cal U_1''=A$, $u^2=1$, $\tau(u)=0$, $\forall u\in \Cal U_1\setminus \{1\}$. Applying Theorem 1.3.2 to $\xi \perp \Cal I_M(A, P)$, 
we get a sequence $\{u_m\}_m \subset \Cal U_1$ such that $\lim_m \|E_P(\xi^* u_m \xi)\|_1=0$. Taking $\alpha>0$ appropriately small 
and $A_0$ to be an $n$-dimensional subalgebra of $A$ generated by projections 
of trace $1/n$ that's $\alpha$-contained in $\{u_m\mid  m_0\leq m\leq m_1\}''$, for some 
$m_0\leq m_1$ sufficiently large so that $\|E_P(\xi^* u_m \xi)\|_1< \alpha$, $\forall m\geq m_0$, one gets condition $3^\circ$ satisfied. 
\hfill$\square$

\heading 3.  Thin factors and MASAs with bounded multiplicity 
\endheading

In the 1950s, W. Ambrose and I.M. Singer have considered 
MASAs in II$_1$ factors $A\subset M$ with the property that the von Neumann algebra $A\vee JAJ \subset \Cal B(L^2M)$, generated by the 
left-right multiplication on $L^2M$ by elements in $A$, is maximal abelian in $\Cal B(L^2M)$. Noticing that this 
is equivalent to $A \vee JAJ$ having a 
cyclic vector (i.e., $\exists \xi \in L^2M$ with $[A\xi A]=L^2M$), this property is analogue to an inclusion of groups $H\subset G$ 
with just one (non-trivial) double co-set over $H$. (Note however that the algebra framework makes it so that 
regular MASAs do satisfy this property (cf. [FM77]), while for a normal subgroup $H\subset G$,  $H \setminus G/H=G/H$ is always large.)

In [Pu61] L. Pukanszky took this  idea further, by noticing that  the type of the algebra $(A \vee JAJ)'\subset \Cal B(L^2M)$ 
is an invariant for the isomorphism class of a MASA inclusion $A\subset M$, and that if an inclusion of groups $H\subset G$ 
with $G$ ICC and $H$ abelian is so that $H\setminus G /H - H$ has $n$ identical classes, 
then the MASA inclusion $A=L(H)\subset L(G)=M$ has the property that $A \vee JAJ$ has multiplicity 
$n$ on $L^2M\ominus L^2A$. In other words, the commutant algebra $(A \vee JAJ)'$, which is always equal to 
$Ae_A\simeq A$ on the reducing space $L^2A$,  is homogeneous of type $\text{\rm I}_n$ 
on  $L^2(M\ominus A)$. Taking appropriate examples $H\subset G$ with $G$ locally finite ICC 
(inspired by a construction in [D54]), he was able to give examples of singular MASAs in the hyperfinite  II$_1$ factor $R$ 
that have ``multiplicity $n$'', for each $1\leq n < \infty$, and are thus mutually non-conjugate by automorphisms of $R$.  

The type of the von Neumann algebra $(A\vee JAJ)' (1-e_A)\subset \Cal B(L^2(M\ominus A))$ (i.e., the list of multiplicities appearing in its 
decomposition as a direct sum of homogeneous type I$_{n_i}$ algebras, $1\leq n_i \leq \infty$) is what one generically calls the {\it Pukanszky invariant}  
of $A\subset M$.  

Of this, we will retain here only the supremum over all the multiplicities $1\leq n_i\leq \infty$ in 
the decomposition $(A\vee JAJ)'=\oplus_i L^\infty(X_i)\otimes M_{n_i\times n_i}(\Bbb C)$. 

\vskip.05in

\noindent
{\bf 3.1. Definitions.} $1^\circ$ Let $M$ be a II$_1$ factor and $A\subset M$ an abelian 
von Neumann subalgebra. We denote by $\text{\rm m}(A\subset M)$ the supremum over all $1\leq m\leq \infty$ with the property that  
$(A \vee JAJ)'$ has a type I$_m$ direct summand, and call it 
the {\it multiplicity of $A\subset M$}. Notice that $\text{\rm m}(A\subset M)= 1$ if and only if $A\vee JAJ$ is maximal abelian in $\Cal B(L^2M)$, 
and that this implies $A$ is maximal abelian  in $M$ (see 3.3 below). An abelian von Neumann subalgebra $A$ in $M$ with the property that 
$A \vee JAJ$ is maximal abelian in $\Cal B(L^2M)$ is called an {\it s-MASA} of $M$.    

$2^\circ$ If $M$ is a II$_1$ factor then we denote $\text{\rm m}_a(M)=\min \{\text{\rm m}(A\subset M)\mid A$ a MASA in $M\}$. 
Thus,  $\text{\rm m}_a(M)=1$ if and only if $M$ has an s-MASA. 
A II$_1$ factor with this property is called an {\it s-thin} factor. 

\vskip.05in 

Let us note right away that the multiplicity of MASAs behaves well to taking tensor products and 
intermediate subfactors: if $A_i\subset M_i$, $i=1, 2$, are inclusions of MASAs, then $\text{\rm m}(A_1\overline{\otimes} A_2 \subset 
M_1\overline{\otimes} M_2)= \text{\rm m}(A_1\subset M_1)\text{\rm m}(A_2\subset M_2)$; also, if $B\subset Q \subset P$ is a MASA 
in $P$ with $Q$ an intermediate factor, then $\text{\rm m}(B\subset Q)\leq \text{\rm m}(B\subset P)$. 

\vskip .05in 

\noindent
{\bf 3.2. Examples.} By  (2.9 in [FM77]), any Cartan MASA in a II$_1$ factor is 
an s-MASA. But by [Pu61], there do exist singular s-MASAs as well. 
For instance, when $M\simeq R$ is the hyperfine II$_1$ factor, 
then besides its (unique by [CFW81]) Cartan subalgebra $D$, 
$R$ contains an s-MASA $A$ that is singular. More precisely, the following example of a singular 
MASA $A\subset R$ from ([D54]) has been shown in ([Pu61]) to be an s-MASA: 
represent $R$ as the group factor 
associated with the amenable ICC group $G$ of 
affine transformation on $\Bbb Q$, with its abelian subgroups $T=\Bbb Q$ (translations), $H=\Bbb Q^*$ (homotheties);   
since $G$ is ICC relative to both $T$ and $H$, $D=L(T), A=L(H)$ are MASAs  in $R=L(G)$;  
since $T$ is normal in $G$, $D$ is a Cartan subalgebra in $R$, while  
since $H$ acts transitively on $T\setminus \{0\}$, 
$A$ follows singular in $R$,  with the vector $\xi=\xi_0 + \xi_1\in \ell^2(G)=L^2R$ cyclic for $A \vee JAJ$, where 
$\xi_0$ is the vector corresponding to the trivial element in $T\subset G$ (so translation by $0$) 
and $\xi_1$ is the element in $T\subset G$ corresponding to translation by $1$.

\vskip .05in

\noindent
{\bf 3.3. Remark.} As we mentioned in 3.1.1$^\circ$ above, 
if $A$ is a  von Neumann subalgebra of $M$, then the condition ``$A\vee J_MAJ_M$  maximal abelian 
in $\Cal B(L^2M)$'' implies that $A$ is a MASA in $M$. To see this, note first that this condition implies $B=A'\cap M$ abelian.  
Indeed, since $e_B \in (A\vee J_MAJ_M)'=A\vee J_MAJ_M$, it follows that $e_B(A\vee J_MAJ_M)e_B$ is maximal abelian 
in $\Cal B(L^2B)$. Since  by the commutativity of $B$ we have $A=J_BAJ_B$, this forces $A=B$. 

The terminology ``s-MASA in $M$'' can thus be viewed as 
emphasizing a {\it strengthening} of  the property of being a MASA in $M$. The prefix ``s'' can also be viewed as hinting to  
the  terminology {\it simple MASA}, which has been sometimes used for abelian von Neumann subalgebras satisfying this property 
(N.B.: this has been the original Ambrose-Singer terminology, 
carried on in [K67], [JP82], [Ge97]). The usage of the adjective ``simple'' for a MASA can however be misleading,  
as (non-trivial) abelian von Neumann algebras  do have non-trivial ideals and are thus not simple as rings... 
The terminology ``simple MASA'' may trigger additional confusion as  it has also 
been used by Takesaki  in [T63], but for a different class of MASAs $A\subset M$,    
via a characterization which  has in fact  been later shown equivalent to $A$ being singular (cf. [H79]).

\vskip .05in 

From this point on, if $S$ is a non-empty subset of a Hilbert space $\Cal H$, 
we will use the notation $[S]$ for the Hilbert subspace generated by $S$, i.e., $[S]=\overline{\text{\rm sp}(S)}$, 
and also for the orthogonal projection of $\Cal H$ onto this space, the difference being always clear from the context. 
Also, if $\xi \in \Cal H$ and  $\Cal H_0 \subset \Cal H$ is a vector subspace, then the notation $\xi\in_\delta \Cal H_0$,  
for some $\delta>0$, means that there exists $\eta\in \Cal H_0$ such that $\|\xi-\eta\|_2<\delta$. 
While if $X\subset \Cal H$, then $X\subset_\delta \Cal H_0$ stands for $\xi\in_\delta \Cal H_0$, $\forall \xi\in X$. 

The following result provides alternative characterizations of MASA-multiplicity.

\proclaim{3.4. Proposition} Let M be a separable $\text{\rm II}_1$ factor, $A\subset M$ an abelian von Neumann subalgebra and $n_0\geq 1$ an integer. 
The following conditions are equivalent: 

\vskip.05in
$1^\circ$ Any type $\text{\rm I}_m$ direct 
summand of $(A \vee JAJ)'$ satisfies $m\leq n_0$. 

$2^\circ$ There exists $X\subset L^2M$ with $|X|\leq n_0$ such that $[A X A]=L^2M$. 

$3^\circ$  $\forall F\subset M$ finite, $\forall \delta>0$, 
there exists $X \subset L^2M$ such that $|X|\leq n_0$ and $F\subset_{\delta} \text{\rm sp} A X A$. 

\vskip.05in

Moreover, if $n_0=1$, then in $2^\circ$ above one can take $X=\{b\}$ with $b=b^*\in M$.  
 
If in addition $A$ is a MASA in $M$, then $1^\circ - 3^\circ$ are also equivalent to: 

\vskip.05in

$4^\circ$ The representation $\text{\rm Ad}(\Cal U(A))\curvearrowright L^2(M\ominus A))$ admits a cyclic 
set of $n_0$ vectors.

\endproclaim

Before proving this result, let us notice the following:

\proclaim{3.5. Lemma} Let $\Cal A\subset \Cal B(\Cal H)$ be an abelian von Neumann algebra acting on the Hilbert space $\Cal H$. 
\vskip.05in
$1^\circ$ The supremum over all $n\leq \infty$ such that $\Cal A'$ has a I$_n$ direct summand is equal to the minimum over all $m\leq \infty $ for which there 
exists $X\subset \Cal H$ with $|X|=m$ and $[\Cal AX]=\Cal H$.

$2^\circ$ For any $\eta_1, \eta_2\in \Cal H$, the set 
$L=\{t\in \Bbb C\setminus \{0\} \mid [\Cal A'(\eta_1+t\eta_2)]\neq [\Cal A'(\eta_1)]\vee [\Cal A'(\eta_2)]\}$ 
is at most countable.
\endproclaim

\noindent
{\it Proof}.  $1^\circ$ This part of the statement is the case ``$B$ abelian'' of the Murray-von Neumann coupling constant theorem (see [vN43]) 
relating a finite von Neumann algebra $B\subset \Cal B(\Cal H)$ with its commutant $B'\subset \Cal B(\Cal H)$, by the 
``factor of multiplicity'' $\text{\rm dim}_B\Cal H$.  

$2^\circ$ This part is just  (Lemma 3.5 in [P82]). 

\hfill $\square$

\vskip.05in
\noindent
{\it Proof of} 3.4. Denote $\Cal A=A \vee JAJ$ and $\Cal B=(A\vee JAJ)'\cap \Cal B(L^2M)$. 
By $3.5.1^\circ$ in the above Lemma, we have $1^\circ \Leftrightarrow 2^\circ$, and we clearly have $2^\circ \Rightarrow 3^\circ$.

Condition $3^\circ$ shows that there exists a sequence of subsets $X_n\subset L^2M$ of cardinality at most $n_0$ 
such that if we denote $p_{n}=[A X_n A]\in \Cal B$ 
then $p_n \rightarrow 1$ in the $so$-topology. By Lemma 3.5.1$^\circ$, each $p_n\Cal Bp_n$ is of finite type 
with all homogeneous type I$_m$ summands satisfying $m\leq n_0$. By ([PP84], or 1.2 in [P94]), this is equivalent to 
having the (probabilistic) Pimsner-Popa index of $\Cal Ap_n \subset p_n \Cal B p_n$ at most equal to $n_0$. Since the definition of this index behaves well 
to limits (see e.g., [PP84]), it follows that the index of $\Cal A\subset \Cal B$ is $\leq n_0$ as well, which in turn means that any type I$_m$ 
direct summand of $\Cal B$ must satisfy $m\leq n_0$. Thus,  $3^\circ \Rightarrow 1^\circ$.

Let us now prove that if there exists $\xi\in L^2M$ such that $[A\xi A]=L^2M$, then there 
exists $b=b^*\in M$ such that $[AbA]=L^2M$. 

Let us first notice that there exists $\eta_0=\eta_0^*\in L^2M$ 
such that $[A\eta_0 A] = [A\xi A]=L^2M$.  Indeed, this follows by noticing that $\Cal A=A \vee JAJ$ satisfies $\Cal A'=\Cal A$ and 
so one can apply Lemma 3.5.2$^\circ$ to 
$\eta_1=\Re \xi$, $\eta_2=\Im \xi$, $\Cal A \subset \Cal B(L^2M)$, to get some $t \in \Bbb R$ such that $\eta_0=\eta_1 + t \eta_2$ 
satisfies $[\Cal A(\eta_0)]=[\Cal A (\eta_1]\vee [\Cal A(\eta_2)]=[\Cal A(\xi)]$. 

Let us also note that if for some $\eta\in L^2M$ and $F\subset L^2M$ a finite set we have $F\subset_\delta \text{\rm sp} A \eta A$, 
then there exists $c=c(\delta; \eta, F)\leq 1$ such that if $\eta'\in L^2M$ satisfies $\|\eta'-\eta\|_2< c$, then we still have $F\subset_\delta  \text{\rm sp} A \eta' A$.

Denote $y_n=e_{[-n, -n+1)}(\eta_0)\eta_0+e_{[n-1,n)}(\eta_0)\eta_0$,  $n\geq 1$, 
and note that $y_n=y_n^*\in M$  
are mutually orthogonal in $L^2M$ with $\|y_n\|\leq 2n$, and that $\eta_0 = \lim_m \Sigma_{k=1}^my_k$ in $L^2M$. 

By letting $b_1=y_1$ and then applying 3.5.2$^\circ$ and the above observation for $m\geq 2$,  
we obtain recursively some scalars $t_m\leq 2^{-m}c(2^{-m-1}, \{y_1, ..., y_m\}, b_{m-1})/2m$ such that $b_m=b_{m-1}+t_my_m$ 
satisfy $\{y_1, ..., y_m\} \subset [A b_m A]$.  Thus, if we denote $b=\lim_m b_m \in (M_h)_1$ (note that this limit exists in operator norm, 
because $\|b_m-b_{m-1}\|\leq 2^{-m}$), then $\|b-b_{m-1}\|_2\leq \|b-b_{m-1}\|\leq c(2^{-m-1}, \{y_1, ..., y_m\}, b_{m-1})$, and thus 
$\{y_1, ..., y_m\}\subset_{2^{-m}} [AbA]$. Thus, all $y_i$ belong to $[AbA]$, and hence so does $\eta_0=\Sigma_i y_i$, implying that $[AbA]\supset [A\eta_0A]=L^2M$.

To prove $4^\circ \Leftrightarrow 2^\circ$ under the condition that $A$ is a MASA in $M$, 
we need to show that that there exists a set $X\subset L^2(M\ominus A)$ with $|X|\leq n_0$, 
such that  $\text{\rm sp} \{u\xi u^* \mid u\in \Cal U(A), \xi\in X\}$ is dense in $L^2(M\ominus A)$. 
But by the ``linearization principle'' (5.1 in [P89]), since $A'\cap M=A$, a set $X\subset L^2(M\ominus A)$ 
is cyclic for $\text{\rm Ad}(\Cal U(A))\curvearrowright L^2(M\ominus A))$ 
if and only $[AXA]=L^2(M\ominus A)$.

\hfill$\square$

\proclaim{3.6. Theorem} Let $M$ be a separable $\text{\rm II}_1$ factor and $n_0\geq 1$ an integer. The following conditions are equivalent: 
\vskip .05in
$1^\circ$ $\text{\rm m}_a(M)\leq n_0$. 

\vskip .05in
$2^\circ$  For any  finite dimensional abelian von Neumann subalgebra
$A_0\subset M$, 
any finite subset $F\subset M$ and any 
$\delta>0$, there exists a finite dimensional abelian von Neumann algebra $A_1\subset M$ containing $A_0$ 
and  $X_1\subset L^2M$ with $|X_1|\leq n_0$ such that $F \subset_{\delta} \text{\rm sp}A_1 X_1 A_1$. 

\vskip .05in
$3^\circ$ There exists a sequence of positive numbers $t_n \searrow 0$ such that each $\text{\rm II}_1$ factor $N=M^{t_n}$ 
satisfies the the following property: 
\vskip .05in 
\noindent 
$(3.6.3)$ For any $F\subset N$ finite and 
$\delta>0$, there exist $A_1\subset N$ finite dimensional abelian 
and $X_1 \subset L^2N$, with $|X_1|\leq n_0$, such that $F \subset_{\delta} \text{\rm sp}A_1 X_1 A_1$.  
\endproclaim 
\noindent
{\it Proof}. We clearly have $1^\circ \Rightarrow 2^\circ$. By applying property $2^\circ$ to 
$A_0=\Bbb Ce + \Bbb C(1-e)$, for projections $e$ in $M$ of trace $\tau(e)=t_n$, we see that $2^\circ \Rightarrow 3^\circ$. 

To prove $2^\circ \Rightarrow 1^\circ$, let $\{x_n\}_n\subset (M)_1$ be a sequence of elements that's $\| \ \|_2$-dense in $(M)_1$. 
We first construct recursively an increasing sequence of finite dimensional abelian von Neumann sbalgebras $A_m\subset M$ 
together with a sequence of subsets $X_m\subset L^2M$, with $|X_m|\leq n_0$, 
such that $\{x_1, ..., x_m\}\subset_{2^{-m}} \text{\rm sp} A_m X_m A_m$. 

Assume $(A_m, X_m)$ have been constructed for $1\leq m \leq n$. Apply $2^\circ$ to 
$F=\{x_1, ..., x_{n+1}\}$, $B_0=A_n$ and $\delta=2^{-n-1}$ 
to get a larger finite dimensional algebra $B_1 \supset B_0$, with a subset $X_{n+1}\subset L^2M$ having at most 
$n_0$ elements, such that if we let $A_{n+1}=B$, then $\{x_1, ..., x_{n+1}\}\subset_{2^{-n-1}} 
\text{\rm sp}A_{n+1}X_{n+1} A_{n+1}$.  

Let now $A=\overline{\cup_n A_n}^w$. By construction, it follows that $\forall F\subset M$ finite and $\varepsilon >0$, there 
exists $X\subset  L^2M$, with $|X|\leq n_0$, 
such that $F\subset_{\varepsilon} \text{\rm sp} A X A$. But then 3.4.3$^\circ$ above implies that $m(A\subset M)\leq n_0$. 

Let us finally prove that $3^\circ \Rightarrow 2^\circ$. Let $F\subset (M)_1$ be a finite set, $B_0\subset M$ a finite dimensional abelian 
von Neumann subalgebra and $\delta>0$. We need to prove that there exists a larger finite dimensional abelian algebra $B_1\supset B_0$ with a 
set $X_1\in L^2M$ of at most $n_0$  vectors such that $F\subset_\delta \text{\rm sp}B_1X_1 B_1$. 
 It is clearly sufficient to prove this for $B_0$ generated by minimal projections $\{q_i\mid 0\leq i \leq m\}$ with $\tau(q_i)=t_n$, 
 $1\leq i \leq m$ and $\tau(q_0)< t_n$, where $n$ is sufficiently large so that 
 $t_n< \delta/100$. Let $v_1, ..., v_m$ be partial isometries in $M$ such that $v_i^*v_i=q_1$ and $v_iv_i^*=q_i$, $1\leq i \leq m$.  

Denote $F'=\{v_i^*yv_j \mid y\in F, 1\leq i,j \leq n\}\subset q_1Mq_1$. 
By property $3^\circ$, there exists a finite dimensional abelian von Neumann algebra $A_1\subset q_1Mq_1$ 
and a subset $X_1\in q_1L^2(M)q_1$ of at most $n_0$ elements, such that 
$F'\subset_{\delta/\sqrt{m}} \text{\rm sp} A_1 X_1 A_1$. For each $\xi_1\in X_1$, denote $\eta_1=\Sigma_{i,j=1}^m v_i\xi_1 v_j^*$. 
Let $Y_1$ be the set of all such $\eta_1$ 
and denote $B_1=\Bbb Cq_0 + \Sigma_{i=1}^m v_iA_1v_i^*$. 
Fix $y\in F$. Since $v_i^*yv_j\in F'$, there exist $a_{ij}\in \text{\rm sp} A_1 X_1 A_1$ such that 
$\|v_i^*yv_j-a_{ij}\|^2_{2,q_1Mq_1} \leq \delta^2/m$. 

Note that, by the definitions of $B_1$ and $\eta_1$, the element $b=\Sigma_{i,j=1}^m v_ia_{ij}v_j^*$ belongs to $\text{\rm sp}B_1 Y_1 B_1$. 
Moreover, by using Pythagoras theorem in $M$, we have  
$$
\|y- b\|^2_2=\|y-\Sigma_{i,j=1}^m v_ia_{ij}v_j^*\|_2^2=\Sigma_{i,j=1}^m \|v_i^*yv_j-a_{ij}\|_2^2
$$
$$
=\Sigma_{i,j=1}^m \|v_i^*yv_j - a_{ij}\|_{2,q_1Mq_1}^2\tau(q_1)\leq m^2\tau(q_1)\delta^2/m = \delta^2\tau(1-q_0).
$$
This shows that $F\subset_\delta \text{\rm sp}B_1 Y_1 B_1$, thus finishing the proof. 
\hfill$\square$

\vskip .05in

Motivated by property $(3.6.3)$ above, we will also consider the following: 

\vskip.05in 
\noindent
{\bf 3.7. Definition.} Let $M$ be a II$_1$ factor. We denote by $\text{\rm wm}_a(M)$ the minimum over all 
cardinalities $1\leq m \leq \infty$ with the property that  
given any finite set $F\subset M$ and any $\varepsilon>0$, there exist a subset $X_1\subset L^2M$ with $|X_1|\leq m$ 
and a finite dimensional abelian $^*$-subalgebra $B_1\subset M$ such that $F\subset_\varepsilon [B_1 X_1 B_1]$. 
The II$_1$ factor $M$ is {\it weak s-thin} if $\text{\rm wm}_a(M)=1$, i.e., if for any finite set $F\subset M$ and any $\varepsilon >0$, 
there exists a finite dimensional abelian von Neumann subalgebra $B_1\subset M$ and a vector $\eta_1\subset L^2M$ such 
that $F\subset_{\varepsilon} \text{\rm sp}B_1\eta_1 B_1$. 

\vskip .05in 

Note that $\text{\rm m}_a(M), \text{\rm wm}_a(M)$ 
are both isomorphism invariants for $M$, that measure the  ``thinness'' of $M$ relative to its 
abelian subalgebras and satisfy $\text{\rm m}_a(M) \geq \text{\rm wm}_a(M)$. Notice also that the invariant $\text{\rm m}_a(M)$ is very much in the spirit of what was called $n$-weak thinness 
in [GP99], which denoted the minimal cardinality $1\leq n \leq \infty$ with the property that there exist hyperfinite von Neumann
subalgebras $R_0, R_1\subset M$ and a set $X\subset L^2M$ with $|X|\leq n$ such that $[R_0 X R_1]=L^2M$. 
More precisely, the minimal such $n$ obviously satisfies $n\leq \text{\rm m}_a(M)$.

\proclaim{3.8. Corollary} $1^\circ$ We have $\text{\rm m}_a(M)=\text{\rm m}_a(M^t)$, for any $t>0$. In particular, if $M$ is an s-thin  factor,   
then its amplifications $M^t$ are s-thin factors, $\forall t>0$. Also, 
if $M$ is weak s-thin, then $M_{n \times n} (M)$ is weak s-thin, $\forall n\geq 1$. 

$2^\circ$ If a $\text{\rm II}_1$ factor $M$ has non-trivial fundamental group $($e.g., if $M$ is a McDuff factor$)$, then $\text{\rm wm}_a(M)=\text{\rm m}_a(M)$. 
In particular, such a $\text{\rm II}_1$ factor is s-thin iff it is weak s-thin. 

$3^\circ$ If a $\text{\rm II}_1$ factor $M$ is generated by an increasing 
sequence of subfactors $M_n\subset M$, then $\text{\rm m}_a(M)\leq \limsup_n \text{\rm m}_a(M_n)$ 
and $\text{\rm wm}_a(M)\leq \limsup_n \text{\rm wm}_a(M_n)$. In particular, if all $M_n$ are  s-thin $($resp. weak s-thin$)$, 
then $M$ is s-thin $($resp, weak s-thin$)$. 
\endproclaim 
\noindent
{\it Proof}. To prove $1^\circ$, note first that if we assume $\text{\rm m}_a(M)\leq n_0$, then by $1^\circ \Rightarrow 3^\circ$ in Proposition 3.6 
there exists $t_n \searrow 0$ such that 
$\text{\rm wm}_a(M^{t_n})\leq n_0$. But then $\text{\rm wm}_a((M^t)^{t_n/t})\leq n_0$, with $t_n/t \searrow 0$, 
and thus by applying $3^\circ \Rightarrow 1^\circ$ in 3.6, 
it follows that $\text{\rm m}_a(M^t)\leq n_0$. 

Part $2^\circ$ is immediate from the equivalence $1^\circ \Leftrightarrow 3^\circ$ in Proposition 3.6. 

Part $3^\circ$ can be easily deduced from the characterizations in Proposition 3.6. To see this, assume first that $\text{\rm m}_a(M_n)\leq n_0$, $\forall n$.  We will directly construct from this 
a MASA $A$ in $M$ such that $\text{\rm m}_a(M)\leq n_0$. 

Let $\{x_k\}_k$ be a sequence 
of elements in the unit ball of $\cup_n M_n$ which is $\| \ \|_2$-dense in $(M)_1$. Assume we have constructed 
finite dimensional 
abelian von Neumann subalgebras $A_1\subset A_2 ... \subset A_m$ in the $^*$-algebra $\cup_n M_n$ together with subsets 
$X_1, ..., X_m\subset \cup_n M_n$, with $|X_i|\leq n_0$, such that 
$\{x_1, ..., x_k\}\subset_{2^{-k}} \text{\rm sp} A_k X_k A_k$, $1\leq k \leq m$. 
Let $K\geq 1$ be large enough such that $\{x_1, ..., x_{m+1}\}\subset M_K$ and 
$A_m \subset M_K$. By applying $3.6.2^\circ$ to $B_0=A_m$, $F=\{x_1, ..., x_{m+1}\}$ and $\delta=2^{-m-1}$, as well 
as the fact that $m_a(M_K)\leq n_0$, 
we get a larger finite dimensional abelian von Neumann subalgebra $A_{m+1}\supset A_m$ in $M_K$ (which can thus be viewed as a subalgebra in $M\supset M_K$), with a set of $\leq n_0$ vectors $X_{m+1} \subset L^2(M_K)\subset L^2M$ with at most $n_0$ elements, 
such that  $\{x_1, ..., x_{m+1}\}\subset_{2^{-m-1}} \text{\rm sp}A_{m+1} X_{m+1} A_{m+1}$. 

Letting now $A=\overline{\cup_m A_m}^w$, it is trivial to see that $A\subset M$ satisfies condition $5^\circ$ in Proposition 3.4, and thus 
satisfies $\text{\rm m}(A\subset M)\leq n_0$, implying that $\text{\rm m}_a(M)\leq n_0$. 

Assume now that $\text{\rm wm}_a(M_n)\leq n_0$, $\forall n$. To prove that $\text{\rm wm}_a(M)\leq n_0$, it is clearly sufficient to show that 3.7 
holds true for any finite subset $F$ in a prescribed $\| \ \|_2$-dense subset $\Cal X \subset (M)_1$. But then the 
inequality is trivial for $F\subset (\cup_n M_n)_1$. 

\hfill$\square$

\heading 4. Singular and semiregular s-MASAs in s-thin factors
 \endheading

We prove in this section that if $M$ has an s-MASA (i.e. $M$ is s-thin), then it has singular and semi-regular 
s-MASAs (in fact, many of them). In other words,  if $M$ admits cyclic MASA actions of $L^\infty([0, 1])$, then it has 
cyclic MASA actions that are relative weak mixing, respectively have a ``large'' relative compact part.  

\proclaim{4.1. Theorem} Let $N$ be a separable s-thin factor and 
$N \hookrightarrow M_n$ be embeddings of $N$ into separable $\text{\rm II}_1$ factors 
such that $N'\cap M_n$ is of type $\text{\rm I}$, $\forall n$. For each $n$, let 
$P_n \subset M_n$ be a von Neumann subalgebra  such that $N\not\prec_{M_n} P_n$.  

$1^\circ$ There exists a singular s-MASA $A\subset N$ such that 
$\Cal N_{M_n}(A)''=A \vee N'\cap M_n$ and $A\not\prec_{M_n} P_n$, $\forall n$.  

$2^\circ$ There exists a semiregular s-MASA $A\subset N$ such that $\Cal N_{M_n}(A)'' \subset N \vee N'\cap M_n$  and 
$A\not\prec_{M_n} P_n$, $\forall n$. 

Moreover, in both $1^\circ$ and $2^\circ$, 
if $M_n'\cap N^\omega\neq \Bbb C$ then one can choose $A$ so that $M_n'\cap A^\omega\neq \Bbb C$ as well. 
\endproclaim 
\noindent
{\it Proof}. We proceed exactly as in the proof of Theorem 2.1, constructing $A$ iteratively, but with an additional 
``local requirement'' which will insure that in the end, besides being singular (resp. semiregular) in $N$, $A$ is an s-MASA in $N$, 
with its Ad-action on $M_n$ being weak mixing relative to $P_n$, $\forall n$. 

Thus, we take $\{e_m\}_m \subset \{ e\in \Cal P(N) \mid \tau(e)\leq 1/2\}$ to be a $\| \ \|_2$-dense sequence. 
Also, we let $\{x_k\}_k\subset (N)_1$ be $\| \ \|_2$-dense in $(N)_1$ and 
for each $M_n$ we choose a sequence $\{x^n_k\}_k\subset (M_n)_1$ that's $\| \ \|_2$-dense in $(M_n)_1$.  

To prove Part 1$^\circ$, we construct recursively 
an increasing sequence of finite dimensional abelian von Neumann subalgebras $A_m\subset N$ together with projections 
$f_m\in \Cal P(A_m)$ and unitary elements $v_m \in \Cal U(A_mf_m)$, $w_m \in \Cal U(A_m)$, as well as a vector $\xi_m \in L^2N$, such that 
for each $1\leq i,j,k \leq m$ we have: 
$$
\|f_m - e_m\|_2 \leq 13 \|e_m - E_{A_{m-1}'\cap N}(e_m)\|_2 \tag 4.1.1
$$
$$
\|E_{A_m'\cap M_k}({x^k_i}^*v_mx^k_j)(1-f_m)\|_2 \leq 2^{-m}, 1\leq i,j,k \leq m \tag 4.1.2
$$
$$
\|E_{A_m'\cap M_k}(x^k_j)-E_{A_m \vee N'\cap M_k}(x^k_j)\|_2 \leq 2^{-m},  \tag 4.1.3
$$
$$
\|E_{P_k}({x^k_i}^* w_m x^k_j)\|_2 \leq 2^{-m}, 1\leq i,j,k \leq m. \tag 4.1.4
$$
$$
\{x_1, ..., x_m\} \subset_{2^{-m}} \text{\rm sp} A_m \xi_m A_m. \tag 4.1.5
$$
 
Assume we have constructed $(A_m, f_m, v_m, w_m, \xi_m)$ satisfying these properties for $m=1, 2, ..., n$. 
By the proof of Theorem 2.1, we can first construct  a finite dimensional abelian algebra $A^1_{n+1}\subset M$ that contains $A_n$, 
with a projection $f_{n+1}\in A^1_{n+1}$ and unitary elements $v_{n+1}\in A^1_{n+1}f_{n+1}$, $w_{n+1}\in A^1_{n+1}$, such that 
conditions $(4.1.1)-(4.1.4)$ are satisfied for $m=n+1$ and with $A^1_{n+1}$ playing the role of $A_{n+1}$. Finally, since $A^1_{n+1}$ 
is contained in the s-thin factor $N$, by 3.6.2$^\circ$ there exists a refinement $A_{n+1}\subset N$ of $A^1_{n+1}$ 
and a vector $\xi_{n+1}\in L^2N$ such that  condition $(4.1.5)$ is satisfied for $m=n+1$. Note that since $A_{n+1}$ contains $A^1_{n+1}$, 
conditions $(4.1.1)-(4.1.4)$ will be satisfied for $m=n+1$. 

If we now denote $A=\overline{\cup_n A_n}^w$, then the same argument as in the proof of Theorem 2.1 shows that 
due to conditions $(4.1.1)-(4.1.4)$ we have $A'\cap M_k=A \vee N'\cap M_k$, 
$A\not\prec_{M_k} P_k$ and $\Cal N_{M_k}(A)=\Cal U(A \vee N'\cap M_k)$, while condition $(4.1.5)$ implies $A\subset N$ satisfies 
condition 3.4.4$^\circ$ with $n_0=1$ and is thus an s-MASA in $N$. 

In turn, to prove part $2^\circ$ we construct recursively 
a sequence of commuting dyadic matrix subalgebras 
$R_m\subset N$ (i.e., $R_m\simeq M_{2^{k_m}\times 2^{k_m}}(\Bbb C)$, for some $k_m\geq 1$), 
with diagonal subalgebras $D_m\subset R_m$, such that if we denote  
$N_m=\vee_{k=1}^m R_k$, $A_m=\vee_{k=1}^m D_k$, 
there exist a projection $f_m$ of trace $1/2$ in $D_m$,  
unitary elements $v_m \in \Cal U(D_mf_m)$, $w_m \in \Cal U(D_m)$, and a vector $\xi_m\in L^2N$, 
such that if we denote $y^k_i=x^k_i-E_{N\vee N'\cap M_k}(x^k_i)$, then the following properties 
are satified for $1\leq i,j,k \leq m$:  
$$
\|E_{A_m'\cap M_k}({y^k_i}^*v_my^k_j)(1-f_m)\|_2 \leq 1/10; \tag 4.1.6
$$
$$
\|f_ny^k_i(1-f_m)\|_2\geq 2\|y^k_i\|_2/5; \tag 4.1.7
$$
$$
\|E_{A_m'\cap M_k}(x^k_j)-E_{A_m \vee N'\cap M_k}(x^k_j)\|_2 \leq 2^{-m};  \tag 4.1.8
$$
$$
\|E_{P_k}({x^k_i}^* w_m x^k_j)\|_2 \leq 2^{-m};  \tag 4.1.9
$$
$$
\{x_1, ..., x_m\} \subset_{2^{-m}} \text{\rm sp} A_m \xi_m A_m. \tag 4.1.10
$$

Assuming we have constructed these objects up to $m=n$, we construct them for $m=n+1$ as follows. 
From the proof of 2.1.2$^\circ$, we can first construct a dyadic matrix algebra $R^1_{n+1}$ that commutes 
with $N_n$, together with a diagonal subalgebra $D^1_{n+1}\subset R^1_{n+1}$, 
a projection $f_{n+1}\in D^1_{n+1}$ of trace $1/2$ and unitaries $v_{n+1}\in D^1_{n+1}f_{n+1}$, $w_{n+1}\in D^1_{n+1}$,  
such that if we denote $A^1_{n+1}=A_n \vee D^1_{n+1}$, 
$N^1_{n+1}=N_n \vee R^1_{n+1}$, then  conditions $(4.1.6)-(4.1.9)$ are satisfied for $m=n+1$, with $A^1_{n+1}\subset N^1_{n+1}$ 
in the role of $A_{n+1}\subset N_{n+1}$.  

Let $\{e_{ij}\}_{i,j\in J}$ be matrix units for $N^1_{n+1}$, 
with $e_{ii}$ generating $A_n$, and denote $F=\{e_{1i}x_ke_{j1} \mid 
1\leq k \leq n+1, i,j\in J\}$. Since $N_0=(N^1_{n+1})'\cap N$ is s-thin (as an amplification of $N$), it follows that there exists an abelian finite 
dimensional von Neumann subalgebra $B_0\subset N_0$ and a vector $\eta_0\in L^2(N_0)$, such that 
$F\subset_\alpha \text{\rm sp}B_0\eta_0 B_0$, where $\alpha=2^{-n-1}/|J|$. Moreover, we may clearly also assume that $B_0$ is dyadic. 
If we now denote $B_1=B_0\vee A_n\in N$ and $\eta_1=\Sigma_{i,j} e_{i1}\eta_0 e_{1j}\in L^2N$, then Pythagoras Theorem implies that 
$\{x_1, ..., x_{n+1}\}\subset_{2^{-n-1}} \text{\rm sp} B_1 \eta_1 B_1$. Finally, we take a dyadic matrix algebra $R^0_{n+1}\subset N_0$ 
having $B_0$ as a diagonal subalgebra and denote $R_{n+1}=R^1_{n+1}\vee R^0_{n+1}$, $D_{n+1}=D^1_{n+1}\vee B_0$, 
$N_{n+1}=N_n \vee R_{n+1}$, $A_{n+1}=B_1=A_n \vee D_{n+1}$, $\xi_{n+1}=\eta_1$.  
It is then immediate to see that all conditions $(4.1.6)-(4.1.10)$ are satisfied for $m=n+1$.

Let $A=\overline{\cup_n A_n}^w=\vee_n D_n$, $R=\overline{\cup_n N_n}^w=\vee_n R_n$. 
Like in the proof of 2.1.2$^\circ$, condition $(4.1.8)$ insures that $A'\cap M_k=A \vee N'\cap M_k$  
while $(4.1.9)$ implies that $A\not\prec_{M_k} P_k$, $\forall k$. 
Also, by the definitions of $A$ and $R$ we see that $\Cal N_{M_k}(A)''\supset R\vee N'\cap M_k$,  
and thus $A$ is semiregular in $N$. On the other hand, 
condition $(4.1.10)$ shows that $\Cal N_{M_k}(A)''\subset N\vee N'\cap M_k$. Finally, condition $(4.1.10)$ combined with 
the case $n_0=1$ of $3.4.4^\circ$ show that 
$A$ is s-MASA in $N$. 

Finally, let us note that in the proof of the existence of singular s-MASAs in $1^\circ$ (resp. of  semi-regular s-MASAs in $2^\circ$), 
if we assume $M_n'\cap N^\omega\neq \Bbb C$, then exactly as in the proof of Theorem 2.1.1$^\circ$ (respectively 2.1.2$^\circ$), one can 
complement the list of conditions in the recursive construction with a condition insuring that $A$ contains non-trivial central sequences of $M_n$. 
\hfill$\square$

\proclaim{4.2. Corollary} Let $M$ be a separable s-thin $\text{\rm II}_1$  factor. 
Then $M$ has uncountably many mutually non-intertwinable singular $($respectively semiregular$)$ s-MASAs. 
Moreover, if $M$ has the property Gamma, then all these MASAs can be taken to contain non-trivial central sequences of $M$.
\endproclaim 
\noindent
{\it Proof}. The argument in the proofs of Corollaries 2.2 and 2.3 works exactly the same way, by using Theorem 4.1 in lieu of Theorem 2.1.  
\hfill$\square$

\heading 5. Final remarks and open problems \endheading

\noindent
{\bf 5.1. Absence of Cartan MASAs versus s-MASAs}. The first  examples  of (separable) II$_1$ factors without Cartan subalgebras 
were obtained by Voiculescu in [Vo96], who used free probability methods to  
prove that the free group factors $L(\Bbb F_n)$ do not have Cartan MASAs. 
It was then realized that a suitable adaptation of the argument in [Vo96]  shows that 
$L(\Bbb F_n)$ doesn't have s-MASAs ([Ge98]), nor even MASAs with finite multiplicity ([GP99]), in fact $\text{\rm m}_a(L(\Bbb F_n))=\infty$. 
Similar arguments can be used to show that any II$_1$ factor of the form $M=N_1 * N_2$,  
with $N_1, N_2$ finitely generated diffuse von Neumann subalgebras of $R^\omega$, satisfies $\text{\rm m}_a(M)=\infty$. 
Indeed \footnote{I 
am  grateful to Dima Shlyakhtenko for pointing out to me this line of arguments.}, this follows 
by combining (4.1 in [GP99]) with 
the lower estimates on free entropy dimension in [Ju01] and the additivity of Voiculescu's free entropy dimension ([Vo96]).

On the other hand, during the last ten years, a large number of results about absence of Cartan MASAs have been 
obtained through deformation rigidity theory ([OP07], [CS11], [CSU11], [PV11], [PV12], [I12]). For instance, it was shown in [PV11] 
that $L(\Bbb F_n)\overline{\otimes} N$ has no Cartan subalgebras for any finite factor $N$. In many ``absence of Cartan MASA'' results 
that are obtained through deformation-rigidity theory,  one actually obtains classes of  II$_1$ factors  $M$ that are {\it strongly solid} in the sense of [OP07], 
i.e.,   the normalizing algebra of any diffuse amenable $B\subset M$ (in particular of any MASA $A\subset M$) is amenable. 
This is notably the case for the free group factors $M=L(\Bbb F_n)$ ([OP07]) and more generally for all factors $L(\Gamma)$ 
arising from non-elementary hyperbolic groups $\Gamma$ ([CS11]).

Notice that no result about automatic amenability of normalizing algebras of MASAs could 
be obtained using free probability,  while absence of s-MASAs 
could not be shown by using deformation rigidity theory! 

Finally, let us point out that absence of Cartan MASAs in a II$_1$ factor $M$ amounts to having no relative compact 
actions by MASAs on $M$, while strong solidity means the relative compact part of any such  
action is amenable. Also,  absence of s-MASAs in $M$ means there are no cyclic actions by MASAs on $M$.

\vskip.05in 
\noindent
{\it 5.1.1. Problem}. It would be interesting to 
find new proofs of absence of s-MASAs in certain factors. 
This is particularly the case for the II$_1$ factor $L(\Bbb F_n)$, where a direct,  ``elementary'' proof seems possible.  

\vskip.05in 
\noindent 
{\it 5.1.2. Problem}. We have no examples of II$_1$ factors with s-MASAs but without Cartan subalgebras. 
One class of factors that may provide such examples are the crossed product factors of the form $M=R \rtimes \Gamma$, 
with $\Gamma=\Gamma_1 \times \Gamma_2$ where $\Gamma_1, \Gamma_2$ are groups in one of the classes 
in [PV11], [PV12], for which one knows that any regular MASA of $M$ is necessarily contained in $R$ (after conjugacy by a unitary). 
Thus, if the $\Gamma$-action on $R$ ``mixes well'' the Cartan MASAs of $R$, then one should be able to show that $R$ cannot contain regular MASAs
of $M$.

\vskip.05in 
\noindent
{\it 5.1.3. Problem}. Is the weak s-thin property equivalent to s-thin ? 
More generally, do we always have $\text{\rm wm}_a(M)=\text{\rm m}_a(M)$? 
We saw that once a factor $M$ has non-trivial fundamental group, then the two properties are 
equivalent, but it is not clear wether this is the case for any II$_1$ factor. 

\vskip.05in 
\noindent
{\it 5.1.4. Problem}. Another question we leave open is whether $\text{\rm m}_a(M)<\infty$ implies 
$\text{\rm m}_a(M)=1$ and whether there are permanence properties relating $\text{\rm m}_a(M)$ with 
the multiplicity invariant $\text{\rm m}_a(N)$ of its subfactors of finite Jones index $N\subset M$. In particular, 
whether $M$ is s-thin if and only if $N$ is s-thin.

\vskip.05in 
\noindent
{\bf 5.2. Local characterization of factors with Cartan MASAs}.  As we have seen above, existence of Cartan MASAs 
is a property of II$_1$ factors, that many II$_1$ factors, such as  
$L(\Bbb F_n)$, do not have. Factors with Cartan MASAs are precisely the ones 
that admit relative compact actions 
of $L^\infty([0, 1])$. Let us call $\Cal C\Cal F$ the class of such II$_1$ factors. If $A$ is a MASA in a II$_1$ factor $M$, 
then there are ways to characterize the regularity property $\Cal N_M(A)''=M$ which does not 
specifically mention the normalizer of $A$. Thus, it is shown in [PS01] that $A$ is Cartan in $M$ 
iff there exists $\Cal U_0=\Cal U_0^*\subset \Cal U(M)$ such that $\overline{\text{\rm sp}\Cal U_0}=M$ 
and $A \ni a \mapsto E_A(uau^*)\in A$ are c.p. maps with discrete (countable) Fubini decomposition, $\forall u\in \Cal U_0$, 
and also iff this is true for $\Cal U_0=\Cal U$. It is also shown in [PV14] that $A$ is Cartan iff $\Cal M=\langle M, A \rangle$ 
has the Kadison-Singer norm-paving property relative to $\Cal A=A\vee JAJ$ and iff there exists a normal 
conditional expectation of $\Cal M$ onto $\Cal A$. 

However, there exists no local, intrinsic characterization of factors $M$ in the class $\Cal C\Cal F$,  that does not 
specifically use the Cartan MASA of $M$. 
Such a characterization would certainly be very interesting. It may be useful in deformation-rigidity theory, 
but also for studying permanence properties of $\Cal C\Cal F$, such as stability to inductive limits, to finite index extension/restriction, or 
to crossed products by amenable groups. 
The criterion $1.4.1 (2)$ may be of help in this direction. A related question is to find an intrinsic, local characterization of factors with unique 
(up to unitary conjugacy) Cartan subalgebra. 

There are reasons to believe that any irreducible subfactor $N\subset M$ of a factor $M$ in the class $\Cal C \Cal F$ has Jones index 
equal to the square norm of a (finite or infinite) bipartite graph. This may even be true for the (possibly larger) class of all s-thin factors. We will 
discuss the motivations behind this conjecture in a future paper.

\vskip.05in 
\noindent
{\bf 5.3. Strengthened singularity}. As shown in ([P81d]), any II$_1$ factor $M$ with the property (T) 
has a MASA $A\subset M$ with the property that the only automorphisms of $M$ that normalize $A$ are  
the automorphism of $M$ implemented by unitaries in $A$. With the terminology in the remark before Corollary 2.2, 
this amounts to $A$ being Aut$(M)$-singular. 
Equivalently, if $\theta: M \simeq N$ is an isomorphism of $M$ onto another II$_1$ factor, 
then $\theta$ is in some sense uniquely determined by its restriction to $A$, $\theta_{|A}$. For this reasons, 
a MASA with this property in a II$_1$ factor $M$ is called {\it super-singular} in $M$ (see 5.1 in [P13]). It is shown in [P13] that, 
besides property (T) factors, the hyperfinite II$_1$ factor 
has super-singular MASAs as well. It is an open problem whether any separable II$_1$ factor has super-singular MASAs. 

The proof of Corollary 2.2 shows that the following property for a MASA $A\subset N$ implies super-singularity: 
given any embedding of $N$ into a II$_1$ factor $M_0$ such that $[M_0: N]<\infty$, 
one has $\Cal N_{M_0}(A)=\Cal U(A) \Cal U(N'\cap M_0)$. By ([P86]), any property (T) II$_1$ factor $N$ has a MASA $A\subset N$ satisfying  
this strengthened super-singularity. Indeed, by [J83] any embedding with finite index $N\hookrightarrow M_0$ arises from a basic construction 
$M_0=\langle N, P \rangle$, for some subfactor $P\subset N$ with $[N: P]=[M_0: N]$, while 
by (Theorem 4.5.1 in [P86]) there are only countably many subfactors of finite index of $N$ up to unitary conjugacy, 
so Theorem 2.1 applies. 

On the other hand, Theorem 2.1 
suggests the following question: does there exist  a separable II$_1$ factor $N$ with a MASA $A\subset N$ having the 
property that for any embedding of $N$ into a separable II$_1$ factor $M_0$ with $N'\cap M_0$ atomic, 
one has $\Cal N_{M_0}(A)=\Cal U(A) \Cal U(N'\cap M_0)$? The answer to this 
question is however negative: if $A\subset N$ is a MASA 
that one also identifies with a Cartan MASA in the hyperfinite II$_1$ factor,  $A\simeq L^\infty([0, 1]) 
\hookrightarrow R$, then $M_0=N*_A R$ 
has the property that $N'\cap M_0=\Bbb C1$ but $\Cal N_{M_0}(A)''\supset R$. 

Another strengthening of the singularity property for a MASA $A\subset M$ is obtained by requiring $A$ to be {\it maximal amenable}  
(or equivalently, maximal AFD, by [C75]) in $M$, i.e., to be so that there exists no intermediate amenable 
subalgebra $A\subset B \subset M$ with $A\neq B$ (so $M$ must be non-amenable). The existence of such 
MASAs was discovered in [P81c], where it was shown that $A=L^\infty([0, 1])$ is maximal amenable in $M=A * P$. In  
particular, the MASA $A_u$ generated by one of the generators of the free group $u\in \Bbb F_n$ is maximal amenable in $M=L(\Bbb F_n)$. 
We have conjectured in the early 1980s that any non-amenable II$_1$ factor contains maximal amenable MASAs. We   
will discuss this problem in details in a forthcoming paper.

\head  References \endhead

\item{[AP17]} C. Anantharaman, S. Popa: ``An introduction to II$_1$ factors'', \newline www.math.ucla.edu/$\sim$popa/Books/IIun-v13.pdf

\item{[CS11]} I. Chifan, T. Sinclair: {\it On the structural theory of} II$_1$ {\it factors of negatively curved groups},  
Ann. Sci. de l'\' Ecole Norm. Sup. {\bf 46} (2013), 1-34.

\item{[CSU12]} I. Chifan, T. Sinclair, B. Udrea: {\it On the structural theory of} II$_1$ {\it factors of negatively curved groups} II, 
Adv. in Math. {\bf 245} (2013), 208-236.

\item{[C75]} A. Connes: {\it Classification of injective factors}, Ann. of Math., {\bf 104} (1976), 73-115.

\item{[CFW81]} A. Connes, J. Feldman, B. Weiss: {\it An amenable equivalence relation is generated by a single
transformation}, Erg. Theory Dyn. Sys.  {\bf 1} (1981), 431-450.

\item{[D54]} J. Dixmier: {\it Sousanneaux abeliens maximaux dans les facteurs de type fini}, Ann. of Math. {\bf 59} (1954), 279-286. 

\item{[D57]} J. Dixmier: ``Les algebres d'operateurs dans l'espace hilbertien'', Gauthier-Vill-\newline ars, Paris 1957, 1969. 

\item{[FM77]} J. Feldman, C.C. Moore: {\it Ergodic equivalence
relations, cohomology, and von Neumann algebras} II, Trans. AMS {\bf 234} (1977), 325-359. 

\item{[G98]} L.M. Ge: {\it Applications of free entropy to finite yon Neumann algebras}, Amer. J. Math. {\bf 119} 
(1997), 467-485.

\item{[GP98]} L.M. Ge, S. Popa:  {\it On some decomposition properties
for factors of type II$_1$}, Duke Math. J., {\bf 94} (1998),
79-101.

\item{[H79]} P. Hahn: {\it Reconstruction of a factor from measures on TakesakiÕs unitary equivalence relation}, 
J. Funct. Anal. {\bf 31} (1979) 263-271.

\item{[J83]} V.F.R. Jones: {\it Index for subfactors}, {\bf 72} (1983), 1-25. 

\item{[JP81]} V.F.R. Jones, S. Popa: {\it Some properties of MASAÕs in factors}, 
in ÒProceedings of VI-th Conference in Operator TheoryÓ, Herculane-Timisoara 1981, 
I. Gohberg (ed.), Birkhauser Verlag, 1982, pp 210-220.

\item{[Ju01]} K. Jung: {\it The free entropy dimension of hyperfinite von Neumann algebras}, 
Trans. AMS {\bf 355} (2003), 5053-5089. 

\item{[K67]} R.V. Kadison: {\it Problems on von Neumann algebras},
Baton Rouge Conference 1967, unpublished notes. 

\item{[MvN36]} F. Murray, J. von Neumann:
{\it On rings of operators}, Ann. Math. {\bf 37} (1936), 116-229.

\item{[MvN43]} F. Murray, J. von Neumann: {\it Rings of operators
IV}, Ann. Math. {\bf 44} (1943), 716-808.

\item{[OP07]} N. Ozawa, S. Popa: {\it On a class of} II$_1$ {\it
factors with at most one Cartan subalgebra}, Annals of Mathematics {\bf 172} (2010),
101-137 (math.OA/0706.3623)

\item{[PP84]} M. Pimsner, S. Popa, {\it Entropy and index for
subfactors}, Annales Scient. Ecole Norm. Sup., {\bf 19} (1986),
57-106.

\item{[P81a]} S. Popa: {\it On a problem of R.V. Kadison on maximal
abelian *-subalgebras in factors}, Invent. Math., {\bf 65} (1981),
269-281.

\item{[P81b]} S. Popa: {\it Orthogonal pairs of *-subalgebras in
finite von Neumann algebras}, J. Operator Theory, {\bf 9} (1983),
253-268.

\item{[P81c]} S. Popa: {\it Maximal injective subalgebras in factors
associated with free groups}, Advances in Math., {\bf 50} (1983),
27-48.

\item{[P81d]} S. Popa:  {\it Singular maximal abelian *-subalgebras in
continuous von Neumann algebras}, J. Funct. Analysis, {\bf 50}
(1983), 151-166. 

\item{[P82]}  S. Popa: {\it Notes on Cartan subalgebras in type
$\text{\rm II}_1$ factors}. Mathematica Scandinavica, {\bf 57} (1985),
171-188.

\item{[P86]} S. Popa: {\it Correspondences}, INCREST Preprint 56/1986, www.math.ucla.edu/ \newline popa/preprints.html

\item{[P89]} S. Popa: {\it Classification of subfactors: the
reduction to commuting squares}, Invent. Math., {\bf 101} (1990),
19-43.

\item{[P92]} S. Popa: {\it Classification of amenable subfactors of
type} II, Acta Mathematica, {\bf 172} (1994), 163-255.

\item{[P94]} S. Popa: {\it Classification of subfactors and of their
endomorphisms}, CBMS Lecture Notes, {\bf 86}, Amer. Math. Soc.
1995.

\item{[P97]} S. Popa: {\it Some properties of the symmetric enveloping algebras
with applications to amenability and property T},
Documenta Mathematica, {\bf 4} (1999), 665-744.

\item{[P01]} S. Popa: {\it On a class of type} II$_1$ {\it factors with Betti numbers invariants}, 
Ann. of Math. {\bf163} (2006), 809-899. 

\item{[P03]} S. Popa: {\it Strong Rigidity of}  II$_1$ {\it Factors
Arising from Malleable Actions of $w$-Rigid Groups} I, Invent. Math.,
{\bf 165} (2006), 369-408. (math.OA/0305306).

\item{[P05a]} S. Popa: ``Deformation-rigidity theory'', NCGOA mini-course, Vanderbilt University, May 2005.  

\item{[P05b]} S. Popa:  {\it Cocycle and orbit equivalence superrigidity
for malleable actions of w-rigid groups}, Invent. Math. {\bf 170}
(2007), 243-295 (math.GR/0512646).

\item{[P06]} S. Popa: {\it Deformation and rigidity for group actions
and von Neumann algebras}, in ``Proceedings of the International
Congress of Mathematicians'' (Madrid 2006), Volume I, EMS Publishing House,
Zurich 2006/2007, pp. 445-479.

\item{[P13]} S. Popa: {\it A} II$_1$ {\it factor approach to the Kadison-Singer problem}, 
Comm. Math. Physics. {\bf 332} (2014), 379-414 (math.OA/1303.1424).

\item{[PS01]} S. Popa, D. Shlyakhtenko: {\it Cartan subalgebras and
bimodule decomposition of type} II$_1$ {\it factors}, Math. Scandinavica
{\bf 92} (2003), 93-102.

\item{[PV11]} S. Popa, S. Vaes:  {\it Unique Cartan decomposition for} II$_1$ 
{\it factors arising from arbitrary actions of free groups}, Acta Mathematica, {\bf 194} (2014), 237-284 

\item{[PV12]} S. Popa, S. Vaes: {\it Unique Cartan decomposition for} II$_1$ {\it factors arising from arbitrary actions of hyperbolic groups}, 
Journal fur die reine und angewandte Mathematik, {\bf 690} (2014), 433-458.

\item{[PV14]} S. Popa, S. Vaes: {\it Paving over arbitary MASAs in von Neumann algebras}, 
Analysis and PDE  (2015) 101-123 (math.OA/1412.0631) 

\item{[Pu61]} L. Pukanszky: {\it On maximal abelian subrings in factors of type} II$_1$,  Canad. J. Math. {\bf 12} (1960), 289-296.

\item{[T63]} M. Takesaki: {\it On the unitary equivalence among the components of decompositions of 
representations of involutive Banach algebras and the associated diagonal algebras}, 
Tohoku Math. J. {bf 15} (1963), 365-393. 

\item{[V07]} S. Vaes: {\it Explicit computations of all finite index bimodules for a family of} II$_1$ {\it factors}, 
Ann. Sci. de l'\' Ecole Norm. Sup.  {\bf 41} (2008), 743-788.

\item{[Vo96]}  D. Voiculescu: {\it The analogues of entropy and Fisher information measure in free probability} III, GAFA {\bf 6} (1996), 172-199. 

\enddocument